\documentclass[a4paper, 12pt, twoside]{article}

\usepackage[a4paper, margin=1in]{geometry}
\usepackage[latin1]{inputenc}
\usepackage[T1]{fontenc}
\usepackage[english]{babel}
\usepackage{textcomp} 
\usepackage{graphicx}
\DeclareGraphicsExtensions{.jpg,.mps,.pdf,.png,.gif}
\usepackage[french]{varioref} 
\usepackage[]{amsmath,amssymb,amsfonts}

\usepackage{float}
\usepackage{bbm}
\usepackage{bbold}
\usepackage{amsthm}
\usepackage{dsfont}
\usepackage{color}
\usepackage{natbib}
\usepackage{pict2e}
\newcommand{\Ben}{\begin{enumerate}}
\newcommand{\Een}{\end{enumerate}}
\newcommand{\Bit}{\begin{itemize}}
\newcommand{\Eit}{\end{itemize}}
\newcommand{\Beq}{\begin{equation}}
\newcommand{\Eeq}{\end{equation}}
\newcommand{\Ba}{\begin{align*}}
\newcommand{\Ea}{\end{align*}}
\newcommand{\Mb}{\mathbf}

\newtheorem{Th}{Theorem}

\newtheorem{Prop}{Proposition}

\newtheorem{Rq}{Remark} 
\newtheorem{Corr}{Corollary}
\newtheorem{Def}{Definition}

\title{Spatial risk measures and applications to max-stable processes}
\begin{document}
\author{
Erwan Koch\footnote{ISFA, CREST and ETH Zurich (Department of Mathematics, RiskLab). erwan.koch@math.ethz.ch}
}

\maketitle

\begin{abstract}
The risk of extreme environmental events is of great importance for both the authorities and the insurance\footnote{In the whole paper, insurance refers in particular to reinsurance.} industry. This paper concerns risk measures in a spatial setting, in order to introduce the spatial features of damages stemming from environmental events into the measure of the risk. We develop a new concept of spatial risk measure, based on the spatially aggregated loss over the region of interest, and propose an adapted set of axioms for these spatial risk measures. These axioms quantify the sensitivity of the risk measure with respect to the space and are especially linked to spatial diversification. The proposed model for the cost underlying our definition of spatial risk measure involves applying a damage function to the environmental variable considered. 
We build and theoretically study concrete examples of spatial risk measures based on the indicator function of max-stable processes exceeding a given threshold. Some interpretations in terms of insurance are provided.

\medskip

\end{abstract}

\section{Introduction}

It is of prime importance for both authorities and insurance companies to take the spatial features of environmental risks into account. For authorities, it is crucial to be able to detect the areas at risk: is it safe to build houses in a given area or better somewhere else? Similarly, an insurance company has to choose its geographical zone of activity as well as its portfolio size. The last issue is obviously related to spatial diversification. Thus, tools (especially risk measures) capable of quantitatively dealing with spatial diversification are needed.

The notion of risk measure has been widely studied in the literature. A risk measure $\Pi$ is a function from a set of random variables (typically a cone) to the real numbers, that satisfies some axioms. The seminal paper by \cite{artzner1999coherent} introduced the concept of coherent risk measure, which was then generalized to the convex case by \cite{follmer2002convex} and \cite{frittelli2002putting}. This static framework for risk measures was then extended to the conditional and the dynamic setting. For a detailed review of conditional and dynamic risk measures, we refer to \cite{acciaio2011dynamic}. The most widely used risk measure in the regulatory context is the Value-at-Risk (VaR).\footnote{Although it would be more rigorous to call this risk measure a quantile, we will mainly use the term VaR since it now entered in the common language of applications in banking and insurance.}

The aforementioned risk measures are univariate. In $\mathds{R}$, the natural order allows to easily define the notion of quantile and therefore VaR. However, in dimensions higher than 1, the lack of such a natural order makes straightforward generalizations non-trivial. This is why many different definitions of multivariate quantiles have emerged in the literature. For a detailed review of these, we refer to \cite{serfling2002quantile}. Regarding extensions of VaR to the multivariate setting, see e.g. \cite{embrechts2006bounds} as well as \cite{cousin2013multivariate}.

To the best of our knowledge, only \cite{follmer2014spatial} and \cite{follmerspatial} use the expression \textit{spatial risk measure}. Consider that a local conditional risk assessment is carried out at each node of a network of financial institutions, in the sense that the risk measure applied takes into account the situation at the other nodes. The main issue they raise is whether the local risk assessments can be aggregated in a consistent way in order to provide a global risk measure.

Since risk measures initially appeared to deal with financial risks, in an insurance context, they do not make explicit the influence of the region where the policies were underwritten. However, in an insurance portfolio, this particular region has an obvious impact on the risk undertaken by the company. In this paper, we introduce a new notion of spatial risk measure by explicitly disentangling the spatial region and the hazard generating losses over this region. Then, we study how the measure of the risk is expected to evolve with respect to some of the features of the spatial region, such as its location and its size. This leads to a set of axioms adapted to the spatial context. Contrary to the axioms proposed by \cite{artzner1999coherent}, we study the sensitivity of the measure of the risk with respect to the space variable.

Let us denote by $A \subset \mathds{R}^2$ the region under consideration and by $\{ C_P(\Mb{x}) \}_{\Mb{x} \in \mathds{R}^2}$ the process of the cost (in the whole paper, cost refers to an economic or insured cost) due to a particular environmental hazard (e.g. a hurricane). The easiest approach to build a spatial risk measure is to integrate a (univariate and static) risk measure (e.g. the variance or VaR) over $A$, i.e. to consider $\dfrac{1}{|A|} \displaystyle \int_A \Pi ( C_P(\Mb{x}) ) \  d\Mb{x},$ where $|.|$ denotes the Lebesgue measure. However, if the process $C_P$ is stationary,\footnote{Throughout the paper, stationarity refers to strict stationarity.} then the distribution of $C_P(\Mb{x})$ is independent of $\Mb{x}$ and thus the previous quantity is equal to $\Pi ( C(\Mb{0}) ).$ The corresponding spatial risk measure reduces to the univariate risk measure associated with a single site, meaning that this approach does not account for the spatial dependence structure of the cost process. In order to overcome this defect, we define our spatial risk measure by applying a (univariate and static) risk measure to the normalized aggregated loss over $A$.

After having defined our general notion of spatial risk measure and the corresponding set of axioms, we introduce a general model for the cost process. The latter involves a mapping of the environmental variable under consideration to a cost via a damage function. Then, the paper focuses on the study of concrete examples of spatial risk measures based on a specific cost process involving the indicator function of threshold exceedances of the environmental variable. In a context of climate change, some extreme events tend to be more and more frequent; see e.g. \cite{SwissRe}. It is of prime importance for authorities as well as for the insurance industry to assess the risk of natural disasters. A precise assessment of the risk of extreme events is crucial in order to satisfy capital requirements under the Solvency II regulatory framework. Therefore, due to the spatial feature of the environmental events, we model the process of the environmental variable using max-stable processes, which constitute an extension of extreme value theory to the level of stochastic processes \citep[see e.g.][]{haan1984spectral, de1986stationary, resnickextreme}.

The remainder of the paper is organized as follows. Our concept of spatial risk measure, its corresponding set of axioms and a general model for the cost process are introduced in Section \ref{Chap_RiskMeasures_Sec_Spatial_Risk_Measures}. Then, Section \ref{Sec_Model_Economic_Loss} describes the specific model for the cost process we use subsequently in order to build concrete examples of spatial risk measures. These examples are studied in Section \ref{Chap_RiskMeasures_Sec_Threshold_Damage_Function}. Section \ref{Chap_RislMeasures_Sec_Conclusion} concludes. Throughout the paper, the elements belonging to $\mathds{R}^d$ for some $d \geq 2$ are denoted using bold symbols, whereas those in more general spaces will be denoted using normal font. All proofs can be found in the Appendix.

\section{Spatial risk measures}
\label{Chap_RiskMeasures_Sec_Spatial_Risk_Measures}

Let $\mathcal{A}$ be the set of all compact subsets of $\mathds{R}^2$ with a positive Lebesgue measure. Denote by $\mathcal{P}$ a family of distributions of real-valued stochastic processes on $\mathds{R}^2$ having locally integrable sample paths. Each process represents the cost caused by the events belonging to specified classes and occurring during a given time period, say $[0,T_L]$. In the following, $T_L$ is considered as fixed and does not appear anymore for the sake of notational parsimony. The events considered here have a spatial extent and thus it is natural to consider a cost process on $\mathds{R}^2$. Each class of events (e.g. a heat wave or a hurricane) will be referred to as a hazard in the following. Let $\mathcal{L}_{\Pi}$ be the set of all real-valued random variables defined on an adequate probability space. A risk measure typically will be some function $\Pi: \mathcal{L}_{\Pi} \mapsto \mathds{R}$.

\subsection{Definitions}

We first give the definition of the normalized spatially aggregated loss, which allows to disentangle the contribution of the space and the contribution of the hazards. Indeed, in the case of an insurance company, the total loss in a portfolio of risks depends on both the region where the policies have been underwritten and the hazards covered in these policies. 

\begin{Def}[Normalized spatially aggregated loss]
\label{Def_normalized_Loss}
For $A \in \mathcal{A}$ and $P \in \mathcal{P}$, the normalized spatially aggregated loss is defined by
\begin{equation}
\label{Chap_RiskMeasure_Eq_Loss_Normalized}
L_N(A,P)= \dfrac{1}{|A|} \displaystyle\int_{A} C_P(\Mb{x}) \  d\Mb{x},
\end{equation}
where the stochastic process $\{ C_P(\Mb{x}) \}_{\Mb{x} \in \mathds{R}^2}$ has distribution $P$.
\end{Def}
The quantity $L(A,P)= \displaystyle \int_{A} C_P(\Mb{x}) \ d\Mb{x}$ is finite since $A$ is compact and $C_P$ has locally integrable sample paths. The random variable $L(A,P)$ corresponds to the total economic or insured loss over region $A$ due to specified hazards and is therefore of interest for spatial risk management. It seems more relevant, for both theoretical study and practical interpretation, to consider the normalized spatially aggregated loss which is a loss per surface unit and can be interpreted in a discrete setting as the loss per insurance policy.

Using the concept introduced in Definition \ref{Def_normalized_Loss}, we now define our notion of spatial risk measure, which makes explicit the contribution of the space in the risk measurement.

\begin{Def}[Spatial risk measure]
\label{Def_Spatial_Risk_Measure}
A spatial risk measure is a function $\mathcal{R}_{\Pi}$ that assigns a real number to any region $A \in \mathcal{A}$ and distribution $P \in \mathcal{P}$:
$$\begin{array}{ccccc}
\mathcal{R}_{\Pi} & : & \mathcal{A} \times \mathcal{P}  & \to & \mathds{R} \\
  &   & (A,P) & \mapsto & \mathcal{R}_{\Pi}(A,P)=\Pi ( L_N(A,P) ),
\end{array}$$
where $L_N(A,P)$ is defined in \eqref{Chap_RiskMeasure_Eq_Loss_Normalized}.
\end{Def}
If the distribution $P$ of the cost process is given,
then the function $\mathcal{R}_{\Pi}(\cdot,P)$ summarizes, for any region belonging to $\mathcal{A}$, the risk caused by the hazards characterized by $P$. In the following, $\mathcal{R}_{\Pi}(\cdot,P)$ will be referred to as the spatial risk measure induced by $P$. For many useful risk measures $\Pi$ (variance, VaR, Expected Shortfall, \dots), this notion of spatial risk measure allows to take (at least) part of the spatial dependence structure of the process $C_P$ into account. However, this is not true in the trivial case of the expectation due to its linearity.

It now appears natural to analyze how $\mathcal{R}_{\Pi}(A,P)$ evolves with respect to $A$ for a given $P$. Some natural properties of $\mathcal{R}_{\Pi}(\cdot,P)$ are described in the set of axioms presented below. The spatial properties of $\mathcal{R}_{\Pi}(\cdot,P)$ depend on both the risk measure $\Pi$ and the probabilistic properties of the cost process characterized by $P$.

\subsection{A set of axioms for spatial risk measures}
\label{Subsec_Axioms_Spatial}

This section provides a set of axioms in the context of the spatial risk measures introduced above. These axioms concern the spatial risk measures properties with respect to the space and not to the cost distribution, the latter being considered as given by the problem at hand.

\medskip

\begin{Def}[Set of axioms for spatial risk measures]
\label{Chapriskmeasures_Def_Axiomatic}
For $A \in \mathcal{A}$, let $\Mb{b}_A$ denote its barycenter. For a fixed $P \in \mathcal{P}$, we define the following axioms for the spatial risk measure induced by $P$:
\medskip
\Ben
\item \textbf{Spatial invariance under translation:} \\ for all $\Mb{v} \in \mathds{R}^2$ and $A \in \mathcal{A},  \ \mathcal{R}_{\Pi}(A+\Mb{v}, P)=\mathcal{R}_{\Pi}(A,P)$, where $A+\Mb{v}$ denotes the region $A$ translated by the vector $\Mb{v}$.
\item \textbf{Spatial sub-additivity:} \\
for all $A_1, A_2 \in \mathcal{A},\  \mathcal{R}_{\Pi}(A_1 \cup A_2, P) \leq \min \{ \mathcal{R}_{\Pi}(A_1, P),\mathcal{R}_{\Pi}(A_2, P) \}$.
\item \textbf{Asymptotic spatial homogeneity of order $\boldsymbol{- \alpha}\boldsymbol{, \alpha \geq 0}$:} \\ for all convex $A \in \mathcal{A}$,
\Beq
\label{Eq_Axiom_Spatial_Homogeneity}
\mathcal{R}_{\Pi}(\lambda A, P) \underset{\lambda \to \infty}{=} K_1+\dfrac{K_2}{\lambda^{\alpha}} + o\left(\frac{1}{\lambda^{\alpha}}\right),
\Eeq
where $\lambda A$ is the area obtained by applying to $A$ a homothety with center $\Mb{b}_A$ and ratio $\lambda >0$, and $K_1 \in \mathds{R}$, $K_2 \in \mathds{R} \backslash \{ 0 \}$. The constants $K_1$ and $K_2$ can depend on $A$.
\Een
\end{Def}
These axioms are well defined. For all $A_1, A_2$ compact, $A_1 \cup A_2$ is compact as a finite union of compact sets. For all $A \in \mathcal{A}$, $\Mb{v} \in \mathds{R}^2$ and $\lambda>0$, $A+\Mb{v}$ and $\lambda A$ are compact as images of compact sets by continuous functions.
 
It is also possible to introduce the axiom of \textbf{spatial anti-monotonicity:}
for all $A_1, A_2 \in \mathcal{A}$, $A_1 \subset A_2 \Rightarrow \mathcal{R}_{\Pi}(A_2, P) \leq \mathcal{R}_{\Pi}(A_1, P)$. It is easily shown that the latter is equivalent to the axiom of spatial sub-additivity.

The axioms proposed in Definition \ref{Chapriskmeasures_Def_Axiomatic} seem relevant when the cost process $C_P$ is stationary. In this case, there is clearly spatial invariance under translation. Our spatial sub-additivity axiom means that the risk associated with the normalized spatially aggregated loss is lower when considering the union of two regions instead of only one of these. It indicates that there is spatial diversification, which appears as a natural property when $C_P$ is stationary. If this axiom is satisfied, an insurance company would be well advised to underwrite policies in both regions $A_1$ and $A_2$. As we will see in the following, the axiom of asymptotic spatial homogeneity\footnote{Connected to the study of the function $\lambda \mapsto \mathcal{R}_{\Pi}(\lambda A,P)$ is, in a temporal context, the analysis of the term structure of risk measures. For related studies, see e.g. \cite{embrechts2005strategic} and \citet{guidolin2006term}.} of order $- \alpha$ seems relevant, especially for $\Pi$ being the expectation, the variance or VaR. 

Although there are some links between our notion of spatial risk measures and financial risk measures as for instance summarized in \cite{Follmer2004}, the inclusion of the space and the $C_p$ process in Definition \ref{Def_Spatial_Risk_Measure} sets our approach rather aside.

For the study of concrete examples of spatial risk measures as well as practical applications, we need to specify a model for the cost process $C_P$. A potential general model is provided in the following section, although the spatial risk measure framework defined above can be used in a much wider setting.

\subsection{A general model for the cost process}

The model we propose for the cost process $\{ C_P(\Mb{x}) \}_{\Mb{x} \in \mathds{R}^2}$ requires two components. The first one concerns the cost generating hazards. We assume that the cost is only due to a unique class of events, i.e. to a unique natural hazard. The latter (e.g. a heat wave or a hurricane) is described by the stochastic process of an environmental variable (e.g. the temperature or the wind speed, respectively), denoted by $\{ Z(\Mb{x}) \}_{\Mb{x} \in \mathds{R}^2}$. We assume that $Z$ is representative of the risk during the whole period $[0, T_L]$.

The second component involves a model mapping the natural hazard into a damage and thus a cost. This model requires both the destruction percentage and the exposure at each location. The destruction percentage is obtained by applying a damage function (also referred to as vulnerability curve in the literature), denoted by $D$, to the natural hazard. This damage function is specific to the type of hazard considered. The exposure process, denoted by $\{ E(\Mb{x}) \}_{\Mb{x} \in \mathds{R}^2}$, can be considered as deterministic and involves especially the demographic, economic and topographic conditions. Finally, the destruction percentage must be multiplied by the exposure, yielding the following model for the cost process:
\Beq
\label{Eq_Economic_Loss_Model}
\left \{ C_P(\Mb{x}) \right \}_{\Mb{x} \in \mathds{R}^2} = \left \{ E(\Mb{x})\  D \left( Z(\Mb{x}) \right) \right \}_{\Mb{x} \in \mathds{R}^2}.
\Eeq
If insured losses are of interest, then the specific terms of the insurance policies should also be taken into account.

Note that many generalizations of the model in \eqref{Eq_Economic_Loss_Model} are possible, such as the introduction of a dependence of $D$ with respect to the location in order to account for the varying type of building in space, the explicit introduction of the temporal aspect and the consideration of several environmental variables (for instance, in case of hurricanes where damages depend on both the wind speed and the rainfall amount, or in case of multi-hazards). 

\begin{Rq}
The presence of $P$ in the right-hand term of \eqref{Eq_Economic_Loss_Model} is implicit: the distribution $P$ of the process $C_P$ indeed depends on the three components of the right-hand term.
\end{Rq}

\section{A specific model for the cost process}
\label{Sec_Model_Economic_Loss}

\subsection{The model}

In order to study concrete examples of spatial risk measures in the following, we consider a particular case of the model for the cost introduced above. For the purpose of this paper, we choose the exposure to be uniformly equal to unity. Moreover, we consider, for a threshold $u>0$, the damage function
$$ 
 D(z) = \mathds{I}_{ \{ z>u \} }, \quad z>0,
$$
which is adapted to the impact of high values of an environmental variable. In order to obtain a more realistic destruction percentage, we could for instance take
$$D(z) = \dfrac{ \mathds{I}_{ \{ z>u \} } }{k_u},$$
where $k_u \geq 1$ depends on $u$. Nevertheless, there is no loss of generality in setting $k_u=1$.

An important focus of our paper is the modeling of losses stemming from extreme events; this is particularly relevant for authorities as well as insurance companies. Hence, for all $\Mb{x} \in \mathds{R}^2$, we consider $Z(\Mb{x})$ to be the maximum of the considered environmental variable at location $\Mb{x}$ over the period of interest $[0, T_L]$. This choice implies that the loss $L(A,P)= \displaystyle \int_{A} C_P(\Mb{x}) \ d\Mb{x}$ corresponds to the spatial aggregation of the costs caused by the worst events happening at each location during the time period considered.
As explained in the next section, in this case, max-stable processes \citep[see e.g.][]{haan1984spectral, de1986stationary, resnickextreme} are ideally suited for modeling purposes. Thus, for the remainder of the paper, we will assume the process $\{ Z(\Mb{x}) \}_{\Mb{x} \in \mathds{R}^2}$ to be max-stable. Moreover, we will often make the classical assumption that $Z$ has standard Fr\'echet margins, i.e., for all $\Mb{x} \in \mathds{R}^2$ and $z>0$,
$$\mathds{P}(Z(\Mb{x}) \leq z)= \exp \left( -\dfrac{1}{z} \right).$$ A max-stable process having standard Fr\'echet margins will be referred to as a simple max-stable process. 

To summarize, the specific cost process we consider in the following is
\Beq
\label{Eq_Specific_Cost_Process}
\left \{ C_P(\Mb{x}) \right \}_{\Mb{x} \in \mathds{R}^2}=\left \{ \mathds{I}_{ \{ Z(\Mb{x})>u \} } \right \}_{\Mb{x} \in \mathds{R}^2},
\Eeq
where $\{ Z(\Mb{x}) \}_{\Mb{x} \in \mathds{R}^2}$ is a max-stable process and $u>0$. Note that $C_P$ has locally integrable sample paths since $| C_P(\Mb{x})| \leq 1$, for all $\Mb{x} \in \mathds{R}^2$.

\subsection{A short introduction to max-stable processes}
\label{Chap_RiskMeasure_Sec_Maxstable}

In the following, ``$\bigvee$'' denotes the supremum and $d \in \mathds{N} \backslash \{0 \}$. Moreover, $\overset{d}{=}$ and $\overset{d}{\to}$ stand for equality and convergence in distribution, respectively.

\subsubsection{Definition and comments}

\begin{Def}[Max-stable process]
\label{Def_Maxstable_Processes}
A stochastic process  $\left \{ G(\Mb{x}) \right \}_{\Mb{x} \in \mathds{R}^d}$ is said to be max-stable if there exist sequences of continuous functions 
$\left( a_T(\Mb{x}), \Mb{x} \in \mathds{R}^d \right)_{T \geq 1}> 0$ and  \\
$\left(  b_T(\Mb{x}), \Mb{x} \in \mathds{R}^d \right)_{T \geq 1} \in \mathds{R}$
such that if $\{ G_t(\Mb{x})\}_{\Mb{x} \in \mathds{R}^d }, t=1, \dots, T,$ are independent replications of $G$,
$$
\left \{ \frac{ \bigvee_{t=1}^T \left \{ G_t(\Mb{x} ) \right \}-b_T(\Mb{x} )}{a_T(\Mb{x} )}  \right \} _{\Mb{x} \in \mathds{R}^d}  \overset{d}{=} \{ G(\Mb{x} ) \}_{\Mb{x} \in \mathds{R}^d}.
$$
\end{Def}

Now, let $\{ T_i(\Mb{x}) \}_{\Mb{x} \in \mathds{R}^d}, i=1, \dots, n,$ be independent replications of a stochastic process $\{ T(\Mb{x}) \}_{\Mb{x} \in \mathds{R}^d}$. Let 
$\left( c_n(\Mb{x}), \Mb{x} \in \mathds{R}^d \right)_{n \geq 1} >0$ and $\left( d_n(\Mb{x}), \Mb{x} \in \mathds{R}^d \right)_{n \geq 1} \in \mathds{R}$ be sequences of continuous functions. If there exists a non-degenerate process $\{ G(\Mb{x}) \}_{\Mb{x} \in \mathds{R}^d}$ such that
\Beq
\label{Eq_Justification_Maxstable}
 \left \{ \frac{\bigvee_{i=1}^n \left \{ T_i(\Mb{x}) \right \} -d_n(\Mb{x})}{c_n(\Mb{x})} \right \}_{\Mb{x} \in \mathds{R}^d} \overset{d}{\rightarrow}  \left \{ G(\Mb{x}) \right \}_{\Mb{x} \in \mathds{R}^d}, \mbox{ for } n \to \infty,
 \Eeq
then $G$ is necessarily max-stable, see e.g. \cite{haan1984spectral}. Therefore, max-stable processes are well suited to model the joint behavior of the temporal maxima at all points in space. 

For practical purposes, the number $n$ of observations over which the maxima are taken depends on the length $T_L$ of the time period considered. We assume that the limit in \eqref{Eq_Justification_Maxstable} has been reached. This assumption of course needs some statistical justification in the examples analyzed. Classically, $T_L=1 \mbox{ year}$ and $n=365$. However, an insurance company can be more interested in the loss due to some particular events than in that corresponding to the worst event in the year. In that case, we can for instance take $T_L=1 \mbox{ week}$. 

\subsubsection{Spectral representations and max-stable models}
\label{Chapriskmeasures_Subsec_Spectralrepresentations}

Different spectral representations of max-stable processes can be found in the literature. Here, we present a very general one which can be found in \cite{strokorb2015tail} and is especially based on \cite{haan1984spectral}, \cite{schlather2002models} and \cite{kabluchko2009spectral}. Let $(\Omega, \mathcal{F}, \nu)$ be a probability space and let $\nu_d$ denote the Lebesgue measure on the Borel $\sigma$-algebra $\mathcal{B}^d$ of $\mathds{R}^d$. Moreover, let $f:\mathds{R}^d \times \Omega \mapsto [0, \infty]$ be a measurable function.
\begin{Th}
\label{Th_General_Spectral_Representation}
Simple max-stable processes that are separable in probability allow for a spectral representation of the form
\Beq
\label{Eq_General_Spectral_Representation}
\{ Z(\Mb{x}) \}_{\Mb{x} \in \mathds{R}^d} \overset{d}{=} \left \{ \bigvee_{i=1}^{\infty} \left \{ U_i V_{\Mb{x}}(S_i) \right \} \right \}_{\Mb{x} \in \mathds{R}^d},
\Eeq
where the $(U_i, S_i)_i$ are the points of a Poisson point process on $(0, \infty) \times E$ with intensity $u^{-2} du \times \mu(ds)$ for some Polish measure space $(E, \mathcal{E}, \mu)$ and the functions $V_{\Mb{x}}: E \mapsto (0, \infty)$ are measurable such that $\int_{E} V_{\Mb{x}}(s) \ \mu(ds)=1$, for each $\Mb{x} \in \mathds{R}^d$. Moreover, any process of the form \eqref{Eq_General_Spectral_Representation} is a simple max-stable process.
\end{Th}
This general spectral representation especially includes two famous classes of representation:
\Ben
\item \textbf{(Mixed) Moving Maxima representation} \\
Let $f$ satisfy
$\displaystyle \int_{\Omega} \int_{\mathds{R}^d} f(\Mb{x}, \omega) \ d\Mb{x} \  \nu(d\omega)=1,$
and let $(E, \mathcal{E}, \mu)=(\mathds{R}^d \times \Omega, \mathcal{B}^d \otimes \mathcal{F}, \nu_d \times \nu)$. Considering the following spectral functions in \eqref{Eq_General_Spectral_Representation},
$$ V_{\Mb{x}} (\Mb{s}, \omega)=f(\Mb{x}-\Mb{s}, \omega), \quad \mbox{for all }(\Mb{s}, \omega) \in \mathds{R}^d \times \Omega \mbox{ and } \Mb{x} \in \mathds{R}^d,$$ yields the so-called Mixed Moving Maxima processes (M3 processes) or Moving Maxima processes (M2 processes) if the random shape function is deterministic (i.e. if $\nu$ charges only one point $\omega_0 \in \Omega$). 
\item \textbf{Stochastic processes-based representation} \\ Let $f$ satisfy 
$\displaystyle \int_{\Omega} f(\Mb{x}, \omega)\  \nu(d\omega)=1, \mbox{ for all } \Mb{x} \in \mathds{R}^d.$
We consider $(E, \mathcal{E}, \mu)=(\Omega,  \mathcal{F}, \nu)$ and $V_{\Mb{x}}(\omega)=f(\Mb{x}, \omega), \mbox{ for all } \omega \in \Omega \mbox{ and } \Mb{x} \in \mathds{R}^d.$ 
It is equivalent to denote
\Beq
\label{Eq_Spectral_Representation_Stochastic_Processes}
V_{\Mb{x}}(.)=V(\Mb{x}),\quad \mbox{for all } \Mb{x} \in \mathds{R}^d,
\Eeq
where $\{ V(\Mb{x}) \}_{\Mb{x} \in \mathds{R}^d}$ is a stochastic process.
\Een
For specific cases included in Theorem \ref{Th_General_Spectral_Representation}, we refer the reader to \cite{de2007extreme}, Corollary 9.4.5 and Theorem 9.6.7, \cite{smith1990max}, and \cite{schlather2002models}, Theorems 1 and 2.

The probabilistic representation \eqref{Eq_General_Spectral_Representation} has led to different parametric models for max-stable processes, some of which are presented in the following. Let $\{ \varepsilon(\Mb{x}) \}_{\Mb{x} \in \mathds{R}^d}$ be a stationary standard Gaussian process with any correlation function $\rho$.

\medskip

\noindent \textbf{The Smith model} \\
\cite{smith1990max} introduces a specific M2 process where the deterministic shape function $f$ is the density of a $d$-variate Gaussian random vector with mean $\Mb{0}$ and covariance matrix $\Sigma$.

\medskip

\noindent \textbf{The Schlather model} \\
\citet{schlather2002models} proposes to set $ V(\Mb{x}) = \sqrt{2 \pi}  \  \varepsilon(\Mb{x})$ in \eqref{Eq_Spectral_Representation_Stochastic_Processes}. All correlation functions stemming from the geostatistical literature can be used, allowing for a rich diversity of behaviors. 

\medskip

\noindent \textbf{The Brown-Resnick model} \\
\citet{kabluchko2009stationary} introduce a generalization of the model introduced by \cite{brown1977extreme}, by taking in \eqref{Eq_Spectral_Representation_Stochastic_Processes} $V(\Mb{x})=\exp \left( W(\Mb{x}) - \dfrac{\sigma_W^2(\Mb{x})}{2} \right)$, where $\{ W(\Mb{x}) \}_{\Mb{x} \in \mathds{R}^d}$ is a Gaussian process with stationary increments and $ \sigma_W^2(\Mb{x}) = \mbox{Var} ( W(\Mb{x}) ),$ for all $\Mb{x} \in \mathds{R}^d$, where Var stands for the variance. This process is referred to as the Brown-Resnick process and its distribution depends only on the semivariogram of $W,$ defined by $\gamma(\Mb{h})=\frac{1}{2} \mbox{Var} ( W(\Mb{x+h})-W(\Mb{x}) )$, $\Mb{h} \in \mathds{R}^d.$
The special case where $W(\Mb{x})=\sigma_{\varepsilon} \  \varepsilon(\Mb{x})$, with $\sigma_{\varepsilon}>0$, leads to the so-called geometric Gaussian process, see e.g. \cite{davison2012statistical}. It should also be noted that the Smith process is a particular case of the Brown-Resnick process; see e.g. \cite{yuen2013crps}, Section 2.2.

\medskip

A correlation function $\rho$ is said to be isotropic if, for all $\Mb{x}_1, \Mb{x}_2 \in \mathds{R}^d$, $\rho(\Mb{x}_1, \Mb{x}_2)$ only depends on $\| \Mb{x}_1-\Mb{x}_2 \|$, where $\|.\|$ denotes the Euclidean distance. In the cases of the Schlather and geometric Gaussian processes, we will compare and contrast the families presented in Table \ref{Table_Correlation_Families}.
\begin{table}[h!]
\begin{center}
\begin{tabular}{|c|c|c|}
\hline
\textbf{Family} & \textbf{Correlation function ($h\geq0$)} & \textbf{Range of validity} \\
\hline 
\hline
Whittle-Mat\'ern & 
$\rho(h)= \frac{2^{1-c_2}}{\Gamma(c_2)} \left( \frac{h}{c_1} \right)^{c_2} K_{c_2}  \left( \frac{h}{c_1} \right)$ & $c_2 > 0$ \\
Cauchy &
$ \rho(h)= \left( 1+ \left( \frac{h}{c_1} \right)^2 \right)^{- c_2}$ & $c_2 > 0$ \\
Powered exponential &
$ \rho(h)= \exp \left( - \left( \frac{h}{c_1} \right)^{c_2} \right)$ & $0 < c_2 < 2$ \\
\hline
\end{tabular}
\caption{Some families of isotropic correlation functions. Here, $c_1$ and $c_2$ are the range and smoothing parameters, $\Gamma$ is the Gamma function and $K_{c_2}$ is the modified Bessel function of the third kind of order $c_2$.}
\label{Table_Correlation_Families}
\end{center}
\end{table}
\noindent The powered exponential correlation function with $c_2=1$ will be referred to as the exponential correlation function.

\subsubsection{Extremal coefficient and mixing properties}
\label{SubSubSec_Extremal_Mixing}

The extremal coefficient \citep[see e.g.][]{schlather2003dependence} is a measure of spatial dependence for max-stable processes. Let $\{ Z(\Mb{x}) \}_{\Mb{x} \in \mathds{R}^d}$ be a simple max-stable process. In the case of two locations, the extremal coefficient function $\Theta$ is defined by
$$
\mathds{P}\left( Z(\Mb{x}_1) \leq u, Z(\Mb{x}_2) \leq u \right) = \exp \left( -\frac{\Theta(\Mb{x}_1, \Mb{x}_2)}{u} \right), \quad \Mb{x}_1, \Mb{x}_2 \in \mathds{R}^d.
$$
If $Z$ is stationary, then $\Theta(\Mb{x}_1, \Mb{x}_2)$ only depends on $\Mb{h}=\Mb{x}_1 - \Mb{x}_2$. Similarly to the correlation function, $\Theta$ is called isotropic if, for all $\Mb{x}_1, \Mb{x}_2 \in \mathds{R}^d$,
$\Theta(\Mb{x}_1, \Mb{x}_2)$ only depends on $\| \Mb{x}_1-\Mb{x}_2 \|$.

Results about spatial diversification we present in Section \ref{Chap_RiskMeasures_Sec_Threshold_Damage_Function} are expressed in terms of the extremal coefficient and it turns out that the latter is directly related to mixing properties of max-stable processes. For all $h \in \mathds{R}$, let us define $r(h)=2-\Theta(h)$, where $\Theta$ is the extremal coefficient function. \cite{kabluchko2010ergodic} prove (Theorems 3.1 and 3.2, respectively) that a simple, stationary and measurable max-stable process $\{ Z(x) \}_{x \in \mathds{R}}$ is mixing if and only if (iff) $\lim_{h \to \infty} r(h)=0$ and is ergodic iff $\lim_{l \to \infty} \dfrac{1}{l} \displaystyle \int_{h=1}^l r(h) \ dh=0$. These results can be extended to $\mathds{R}^d$; see e.g. \cite{DombryHDR2012}, p.20.

Let us now consider the process $\{ H(\Mb{x}) \}_{\Mb{x} \in \mathds{R}^d}=\{ D ( Z(\Mb{x}) ) \}_{\Mb{x} \in \mathds{R}^d}$, where $D$ is a real-valued and measurable function defined (at least) on the state space of $Z$. It is worth mentioning that if $Z$ is mixing, this also holds true for $H$. Indeed, the $\sigma$-algebra generated by $H$ is included in that generated by $Z$. This means that the mixing properties of the environmental process are also valid for the cost process $C_P$.

\section{Some concrete examples of spatial risk measures}
\label{Chap_RiskMeasures_Sec_Threshold_Damage_Function}

In this section, we consider the normalized spatially aggregated loss obtained by combining \eqref{Chap_RiskMeasure_Eq_Loss_Normalized} and \eqref{Eq_Specific_Cost_Process}:
\Beq
\label{Chapriskmeasures_Excursion_Sets}
L_N(A,P) = \frac{1}{|A|} \int_{A} \mathds{I}_{ \{ Z(\Mb{x}) > u \} } \ d\Mb{x},
\Eeq
where $\{ Z(\Mb{x}) \}_{\Mb{x} \in \mathds{R}^2}$ is a max-stable process and $u>0$. Note that, for all $A \in \mathcal{A}$ and $P \in \mathcal{P}$, $|L_N(A,P)| \leq 1$, meaning that all moments of $L_N(A,P)$ are finite. This quantity is interesting when analyzing, for instance, the impact of high temperatures on populations as well as on the distribution network of electricity, typically as in the case of the European heatwave in 2003. If u is well chosen, $L_N(A,P)$ represents indeed the proportion of the surface area at which the temperature exceeds a dangerous threshold for populations and electric cables.

It should be noted that $\displaystyle \int_{A} \mathds{I}_{ \{ Z(\Mb{x}) > u \} } \ d\Mb{x}$  corresponds to the so-called intrinsic volume of the excursion set $E_u(Z,A)= \{ \Mb{x} \in A : Z(\Mb{x}) \geq u \}$. Excursion sets of stochastic processes have been widely studied in the literature; see e.g. \cite{adler2013high}, \cite{spodarev2014limit} and references therein.

The dependence of $L_N(A,P)$ with respect to $P$ lies in the distribution of $Z$. In the following, different max-stable models, especially those introduced in Section \ref{Chapriskmeasures_Subsec_Spectralrepresentations}, will be considered. In each case, the model used will be explicitly indicated and therefore, we make the dependence in $P$ implicit in $L_N(A,P)$: from now on, $L_N(A,P)$ will be denoted $L_N(A)$. Likewise, the spatial risk measure $\mathcal{R}_{\Pi}(A,P)$ will be denoted $\mathcal{R}_{\Pi}(A)$.

The case of the expectation $\mathcal{R}_1(A) = \mathds{E} \left( L_N(A) \right)$ is trivial if $Z$ has identical margins. Using the linearity property of the expectation, we immediately show, if $Z$ is simple, that, for all $A \in \mathcal{A}$,
$$\mathcal{R}_1(A) = 1 - \exp \left( -\dfrac{1}{u} \right),$$
meaning that $\mathcal{R}_1$ does not depend on the region considered (and thus on its size). The spatial risk measure $\mathcal{R}_1$ obviously satisfies the axioms presented in Definition \ref{Chapriskmeasures_Def_Axiomatic}, the order of the asymptotic spatial homogeneity being 0.\footnote{This is true for all cost processes having identical margins and finite expectation.} The expectation is not a very useful risk measure since it does not involve any information relative to the variability. Furthermore, due to linearity, it does not account for the spatial dependence of the cost process $C_P$. In the following, we study the case of the variance in details, before providing some insights concerning VaR and briefly discussing some practical considerations.

\subsection{The variance}

We consider in this section the quantity $\mathcal{R}_2(A)=\mbox{Var}\left( L_N(A) \right).$ The variance allows taking into account part of the spatial dependence in the risk assessment. Hence, its study is interesting for both the risk management of extreme spatial events and the understanding of some properties of max-stable processes. As we will see, $\mathcal{R}_2$ is linked with the notions of spatial diversification, ergodicity and mixing.

In the following, our aim is to study whether $\mathcal{R}_2$ satisfies the axioms presented in Definition \ref{Chapriskmeasures_Def_Axiomatic}. Before doing so, we study the function $\lambda \mapsto \mathcal{R}_2(\lambda A)$ in detail. It is of course related to the property of asymptotic spatial homogeneity. However, we are also interested in the behavior of $\mathcal{R}_2(\lambda A)$ for finite values of $\lambda$. The latter can be of practical relevance for the insurance industry and, moreover, the results obtained will be used to prove the axiom of spatial sub-additivity. 

\subsubsection{Study of $\mathcal{R}_2(\lambda A), \lambda >0$}

The expression of $\mathcal{R}_2(\lambda A)$ in the case of simple max-stable processes is given in the next theorem. 
\begin{Th}
\label{Prop_Stationary_Risk_Measure_Generalcase}
Let $\left \{ Z(\Mb{x}) \right \}_{\Mb{x} \in \mathds{R}^2}$ be a simple max-stable process having $\Theta$ as extremal coefficient function. For all $A \in \mathcal{A}$ and $\lambda>0$, we have
\Beq
\label{Chapriskmeasures_Eq_Stationary_Risk_Measure_Generalcase}
\mathcal{R}_2(\lambda A)= \frac{1}{\lambda^4 |A|^2} \int_{\lambda A}  \int_{\lambda A} \left[ \exp\left( -\frac{\Theta(\Mb{x}, \Mb{y})}{u} \right) - \exp \left( -\frac{2}{u} \right) \right] \  d\Mb{x}\  d\Mb{y}.
\Eeq
\end{Th}
From Theorem \ref{Prop_Stationary_Risk_Measure_Generalcase}, we can derive the expression of $\mathcal{R}_2(\lambda A)$ in the case of simple max-stable processes having an isotropic extremal coefficient function, when $A$ is either a disk or a square. Note that in the whole paper, disk and square refer to a closed disk with positive radius and a (closed) square with positive side, respectively. The result is given in the next corollary.

\begin{Corr}
\label{Prop_Stationary_Risk_Measure}
Let $\left \{ Z(\Mb{x}) \right \}_{\Mb{x} \in \mathds{R}^2}$ be a simple max-stable process having an isotropic extremal coefficient function $\Theta(h), h \geq 0$.
Then:
\Ben
\item Let $A$ be a disk with radius $R$. For all $\lambda >0$, we have
\Beq
\label{Risk_Measure_Disk}
\mathcal{R}_2 (\lambda A) = - \exp \left( -\frac{2}{u} \right) + \int_{0}^{2R} f_d(h,R) \  \exp \left( -\frac{\Theta(\lambda h)}{u} \right) \  dh,
\Eeq
where $f_d$ is the density of the Euclidean distance between two points independently and uniformly distributed on $A$, given by
$$
f_d(h,R) = \frac{2h}{R^2} \left( \frac{2}{\pi} \arccos \left( \frac{h}{2 R} \right) - \frac{h}{\pi R} \sqrt{1-\frac{h^2}{4 R^2}}   \right).
$$
\item Let $A$ be a square with side $R$. For all $\lambda >0$, we have
\Beq
\label{Risk_Measure_Square}
\mathcal{R}_2 (\lambda A) = - \exp \left( -\frac{2}{u} \right) + \int_{0}^{\sqrt{2} R} f_s(h,R) \ \exp \left( -\frac{\Theta(\lambda h)}{u} \right) \  dh,
\Eeq
where $f_s$ is the density of the Euclidean distance between two points independently and uniformly distributed on $A$, given by: \\
\Bit
\item For $h \in [0,R]$,
$$ f_s(h,R)=\frac{2 \pi h}{R^2} - \frac{8 h^2}{R^3} + \frac{2 h^3}{R^4}.$$
\item For $h \in  [R, R \sqrt{2}]$,
\begin{align*}
& \quad \ f_s(h,R) \\&= \left( -2-b+ 3 \sqrt{b-1} + \frac{b+1}{\sqrt{b-1}} +2 \arcsin \left( \frac{2-b}{b} \right) - \frac{4}{b \sqrt{1-\frac{(2-b)^2}{b^2}}}\right) \frac{2h}{R^2},
\end{align*}
where $b=\dfrac{h^2}{R^2}$.
\Eit
\item In both cases, if $\lim_{\lambda \to \infty} \Theta(\lambda)$ exists, then $\mathcal{R}_2 (\lambda A)$ converges as $\lambda \to \infty$ to the limiting risk measure given by
\Beq
\label{Stationary_Risk_Measure}
- \exp \left( -\frac{2}{u} \right) + \exp \left( -\frac{\lim_{\lambda \to \infty} \Theta(\lambda)}{u} \right).
\Eeq
\Een
\end{Corr}
\begin{Rq}
\label{Chapriskmeasures_Rq_Shape}
For more complex geometric shapes of the region A, little is known on the density of the distance between two points independently and uniformly distributed on $A$ \citep[see e.g.][Section 4.3.3]{moltchanov2012distance}, making this study more difficult.
\end{Rq}
The following result directly follows from Corollary \ref{Prop_Stationary_Risk_Measure}.
\begin{Corr}
Let $Z$ be as in Corollary \ref{Prop_Stationary_Risk_Measure} and $A$ be a disk or a square. If, for all $h \geq 0, \ \Theta(h)=1$ (case of asymptotic perfect dependence whatever the distance), we have
$$
\mathcal{R}_2(A) = \exp \left( -\frac{1}{u} \right) - \exp \left( -\frac{2}{u} \right).
$$
\end{Corr}

Let us consider the max-stable models introduced in Section \ref{Chapriskmeasures_Subsec_Spectralrepresentations}, in the case where they have an isotropic extremal coefficient function. This requires the correlation function $\rho$ and the semivariogram $\gamma$ to be isotropic. Moreover, let us assume that $\rho$ and $\gamma$ are non-increasing and non-decreasing, respectively. In this case, for $h > 0$, the function $\lambda \mapsto \Theta(\lambda h)$ is non-decreasing (see below), giving that $\lambda \mapsto \mathcal{R}_2(\lambda A)$ is non-increasing for $A$ being a disk or a square. Consequently, there is spatial diversification and $\mathcal{R}_2(\lambda A)$ has a limit as $\lambda \to \infty$. Corollary \ref{Prop_Stationary_Risk_Measure} may be of interest for the insurance industry since it allows determining the characteristic dimension of the geographical area required to reach a specified low variance level. As shown in the following corollary, in some cases, the diversification can be total, in the sense that $\lim_{\lambda \to \infty} \mathcal{R}_2(\lambda A) = 0$. 

\begin{Corr}
\label{Chapriskmeasures_Corr_Asymptotic_Independence}
Let $Z$ be as in Corollary \ref{Prop_Stationary_Risk_Measure} and $A$ be a disk or a square. If $\lim_{\lambda \to \infty} \Theta(\lambda)=2$ (case of asymptotic independence when the distance tends towards infinity), we have
\Beq
\label{Eq_Convergence_Variance_Corr_3}
\lim_{\lambda \to \infty} \mathcal{R}_2(\lambda A) = 0.
\Eeq
Moreover, if $\lim_{\lambda \to \infty} \Theta(\lambda)$ exists, then $\lim_{\lambda \to \infty} \Theta(\lambda)=2$ iff $\lim_{\lambda \to \infty} \mathcal{R}_2(\lambda A) = 0$.
\end{Corr}
This result is not surprising. Indeed, as mentioned in Section \ref{SubSubSec_Extremal_Mixing}, if additionally $Z$ is measurable and stationary, $\lim_{\lambda \to \infty} \Theta(\lambda)=2$ implies that $Z$ is mixing. Thus, the process $\mathds{I}_{ \{ Z>u \} }$ is also mixing and therefore mean-ergodic, which implies \eqref{Eq_Convergence_Variance_Corr_3} for all $A \in \mathcal{A}$. Furthermore, from Corollary \ref{Chapriskmeasures_Corr_Asymptotic_Independence}, we deduce, as can be done from \cite{kabluchko2010ergodic} for max-stable processes on $\mathds{R}$, that in the case of a simple, stationary and measurable max-stable process having an isotropic extremal coefficient function, if $\lim_{\lambda \to \infty} \Theta(\lambda)$ exists, then ergodicity and mixing are equivalent. Indeed, if $Z$ is ergodic, then $\mathds{I}_{ \{ Z>u \} }$ is also ergodic and thus mean-ergodic, which implies that $\lim_{\lambda \to \infty} \mathcal{R}_2(\lambda A)=0$. Thus, if $\lim_{\lambda \to \infty} \Theta(\lambda)$ exists, this gives that $\lim_{\lambda \to \infty} \Theta(\lambda)=2$, i.e. that $Z$ is mixing. 

Theorem \ref{Prop_Stationary_Risk_Measure_Generalcase} shows that the behavior of the function $\lambda \mapsto \mathcal{R}_2(\lambda A)$ is mainly driven by the extremal coefficient function $\Theta$; the latter naturally depends on the max-stable model under consideration. In the following, we study the influence of $\lambda$ for some classical simple max-stable processes having an isotropic extremal coefficient function. Since the integrals in Corollary \ref{Prop_Stationary_Risk_Measure} have no closed form, we use a Riemann approximation. Without loss of generality, we set $R=1$. The choice of the threshold $u$ has no influence on the shape of the function $\lambda \mapsto \mathcal{R}_2(\lambda A)$ and we choose $u=1$.

\medskip

\noindent \textbf{The Smith model}

\medskip
 
\begin{figure}[h!]
\center
\includegraphics[scale=0.52]{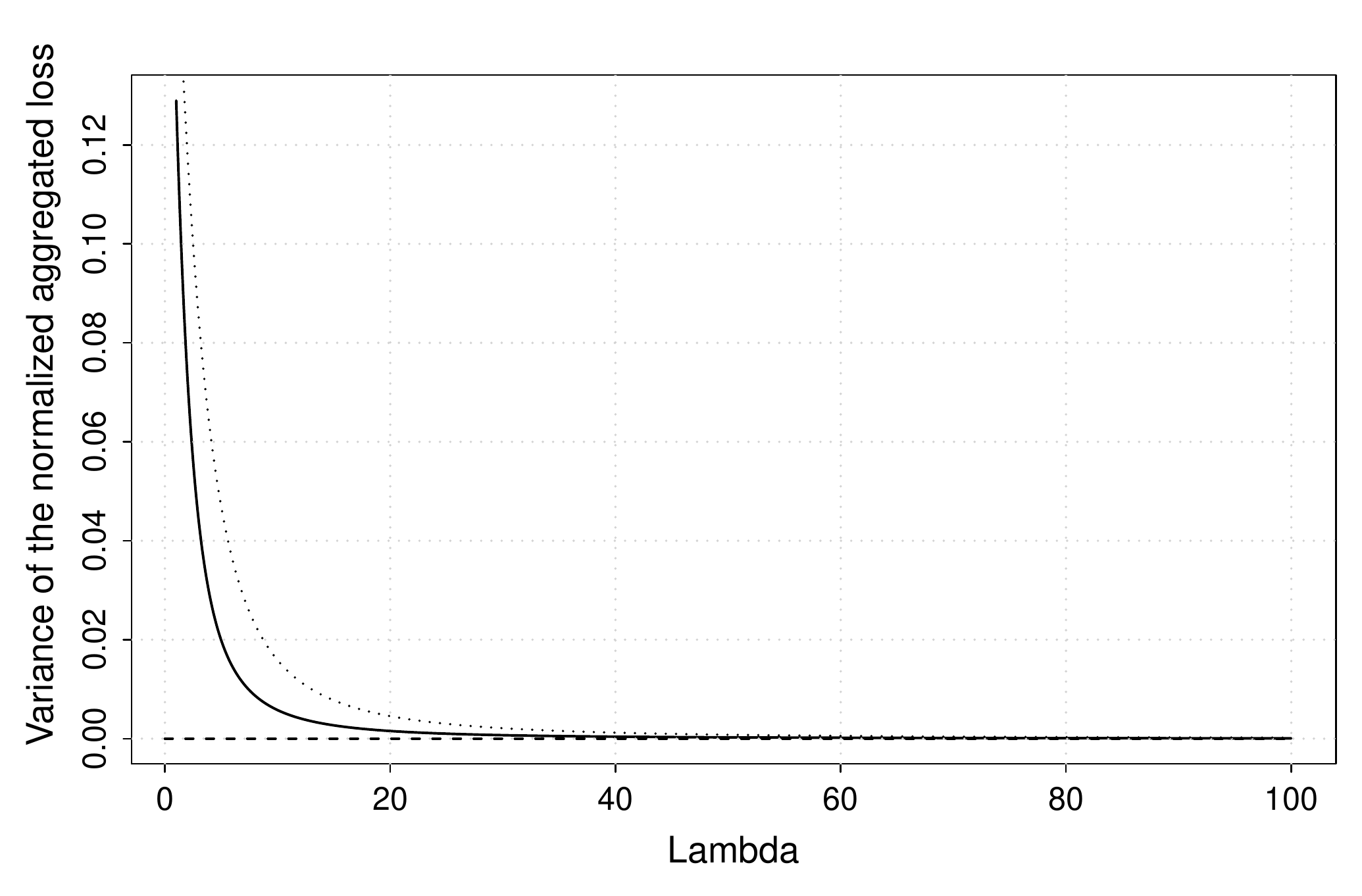}
\caption{The solid (respectively dotted) line depicts the evolution of $\mathcal{R}_2(\lambda A)$ with respect to $\lambda$ in the case of the Smith model, where $A$ is a disk (respectively a square). The dashed line represents the limiting risk measure.}
\label{R2 in case of Smith model}
\end{figure}

The extremal coefficient function is given by \citep[see e.g.][]{davison2012statistical}
$$ \Theta(\Mb{x}_1, \Mb{x}_2)= 2 \Phi \left( \frac{\sqrt{(\Mb{x}_1-\Mb{x}_2)^{'} \Sigma^{-1} (\Mb{x}_1-\Mb{x}_2) }}{2}  \right), \quad \Mb{x}_1, \Mb{x}_2 \in \mathds{R}^2,$$
where $\Phi$ denotes the distribution function of a standard Gaussian random variable and the prime means transposition. The function $\Theta$ is isotropic iff $\Sigma$ is proportional to the identity matrix. Without loss of generality, let $\Sigma$ be the identity matrix.
In this case, we have
$\Theta(h) = 2 \Phi\left( \frac{h}{2} \right)$, for all $h \geq 0$.
Hence, $\lim_{\lambda \to \infty} \Theta(\lambda)=2$
and Corollary \ref{Chapriskmeasures_Corr_Asymptotic_Independence} gives that $\lim_{\lambda \to \infty} \mathcal{R}_2(\lambda A)=0$, meaning that the spatial diversification is total.
We observe in Figure \ref{R2 in case of Smith model} that $\mathcal{R}_2(\lambda A)$ rapidly decreases to the limiting risk measure when $\lambda$ increases.

\medskip

\noindent \textbf{The Schlather model}

\medskip

The extremal coefficient function is given by \citep[see e.g.][]{davison2012statistical}
$$\Theta(\Mb{x}_1, \Mb{x}_2)=1+\sqrt{\dfrac{1-\rho(\Mb{x}_1-\Mb{x}_2)}{2}}, \quad \Mb{x}_1, \Mb{x}_2 \in \mathds{R}^2.$$
The function $\Theta$ is isotropic iff the correlation function $\rho$ is isotropic. In this case, $\Theta(h)= 1+\sqrt{\frac{1-\rho(h)}{2}}$, for all $h \geq 0$. Moreover, for all $h \geq 0$, $\rho(h) \geq -\frac{1}{2}$ \citep[see e.g.][Theorem 3.10]{abrahamsen1997review}, implying that, for all $\lambda \geq 0$, $\Theta(\lambda) < 1.87$, and thus $\lim_{\lambda \to \infty} \mathcal{R}_2(\lambda A)>0$. This shows that the process $\mathds{I}_{ \{ Z>u \} }$ is not mean-ergodic, which is consistent with the fact that the Schlather process is neither mixing nor ergodic. In terms of insurance, this result means that the spatial diversification is never total. In the case of the correlation functions in Table \ref{Table_Correlation_Families}, we have $\lim_{\lambda \to \infty} \rho(\lambda) = 0$, giving, using \eqref{Stationary_Risk_Measure}, that
$\lim_{\lambda \to \infty} \mathcal{R}_2(\lambda A) =  - \exp \left( -\frac{2}{u} \right) +  \exp \left( -\frac{1+\sqrt{\frac{1}{2}}}{u} \right).$

We set the range parameter $c_1=1$ and the smoothing parameter $c_2=0.5$.
\begin{figure}[h!]
\center
\includegraphics[scale=0.56]{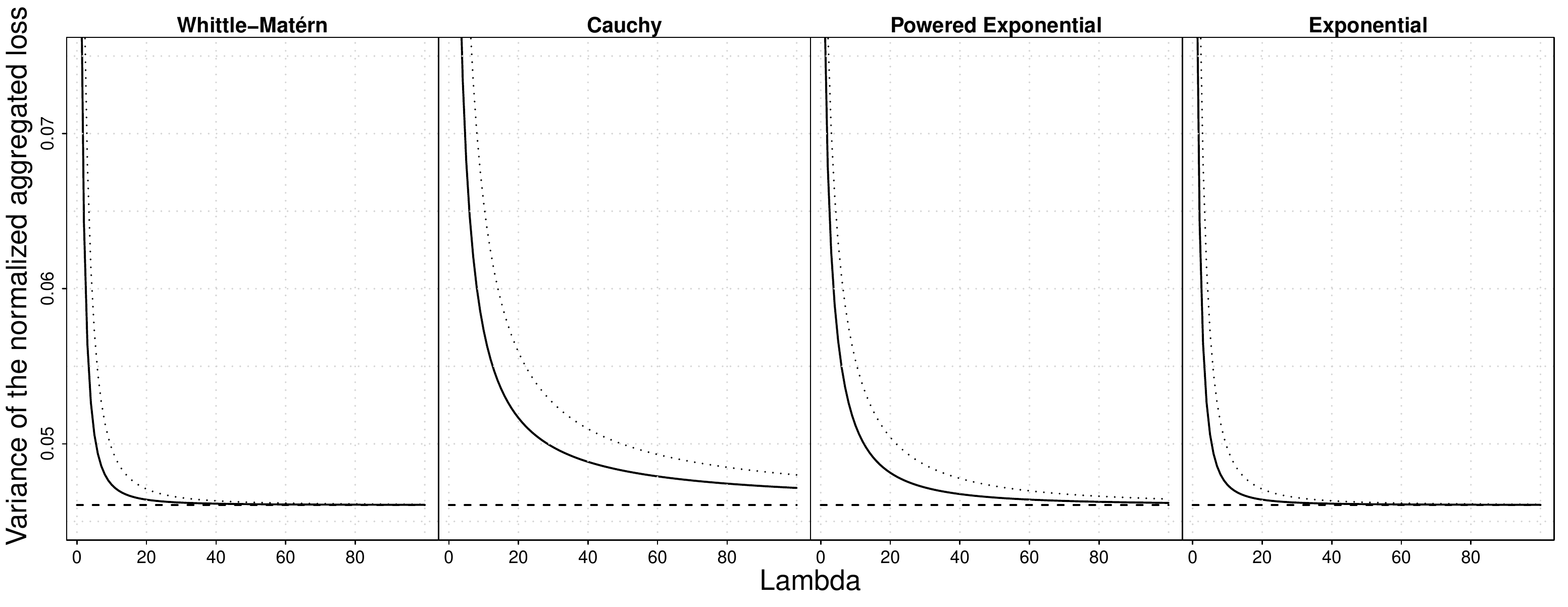}
\caption{Each panel corresponds to a different correlation function. The solid (respectively dotted) line depicts the evolution of $\mathcal{R}_2(\lambda A)$ with respect to $\lambda$ in the case of the Schlather model, where $A$ is a disk (respectively a square). The dashed line represents the limiting risk measure.}
\label{R2 in case of Schlather model for $c_2=1$}
\end{figure}
Figure \ref{R2 in case of Schlather model for $c_2=1$} shows that the speed of decrease of $\mathcal{R}_2(\lambda A)$ to the limiting risk measure depends on the type of correlation function. This decrease is much slower in the cases of the Cauchy and the powered exponential functions. The smoothing parameter also plays an important role, as can be seen from the comparison between the powered exponential case and the exponential one. Hence, the Schlather process allows for a large variety of spatial diversification behavior.

On the whole, we observe that the decrease to the limiting risk measure is slower than in the case of the Smith model. This comparison must be made with equivalent characteristic distances of spatial dependence, meaning that the eigenvalues of $\Sigma$ must be close to $c_1$, which is the case. Note that the decrease is obviously slower when increasing the range parameter $c_1$.

\medskip

\noindent \textbf{The geometric Gaussian model} 

\medskip

The extremal coefficient function is given by \citep[see e.g.][]{davison2012statistical}
 $$ \Theta(\Mb{x}_1, \Mb{x}_2)= 2 \Phi \left( \sqrt{\frac{\sigma_{\varepsilon}^2 [1-\rho(\Mb{x}_1-\Mb{x}_2)]}{2}} \right), \quad \Mb{x}_1, \Mb{x}_2 \in \mathds{R}^2.$$
The function $\Theta$ is isotropic iff $\rho$ is isotropic. In this case, $\Theta(h)= 2 \Phi \left( \sqrt{\frac{\sigma_{\varepsilon}^2 [1-\rho(h)]}{2}} \right)$, for all $h \geq 0$.
Therefore, as before, $\lim_{\lambda \to \infty} \mathcal{R}_2(\lambda A)>0$. In the cases considered here, $\lim_{\lambda \to \infty} \rho(\lambda) = 0$, so that by \eqref{Stationary_Risk_Measure},
$\lim_{\lambda \to \infty} \mathcal{R}_2(\lambda A) =  - \exp \left( -\frac{2}{u} \right) +  \exp \left( -\frac{2 \Phi \left( \sqrt{\frac{\sigma_{\varepsilon}^2}{2}} \right) }{u} \right).$

We set $\sigma_{\varepsilon}=1$ and, as previously, $c_1=1$ and $c_2=0.5$.
From the plots (available upon request), we can draw very similar conclusions to those related to the Schlather model. The main difference consists in the value of the limiting risk measure.

In order to allow for faster spatial diversification, we consider a new instance of a M2 process, the so-called tube model, defined below.

\medskip

\noindent \textbf{The tube model} 

\medskip

\begin{figure}[h!]
\center
\includegraphics[scale=0.52]{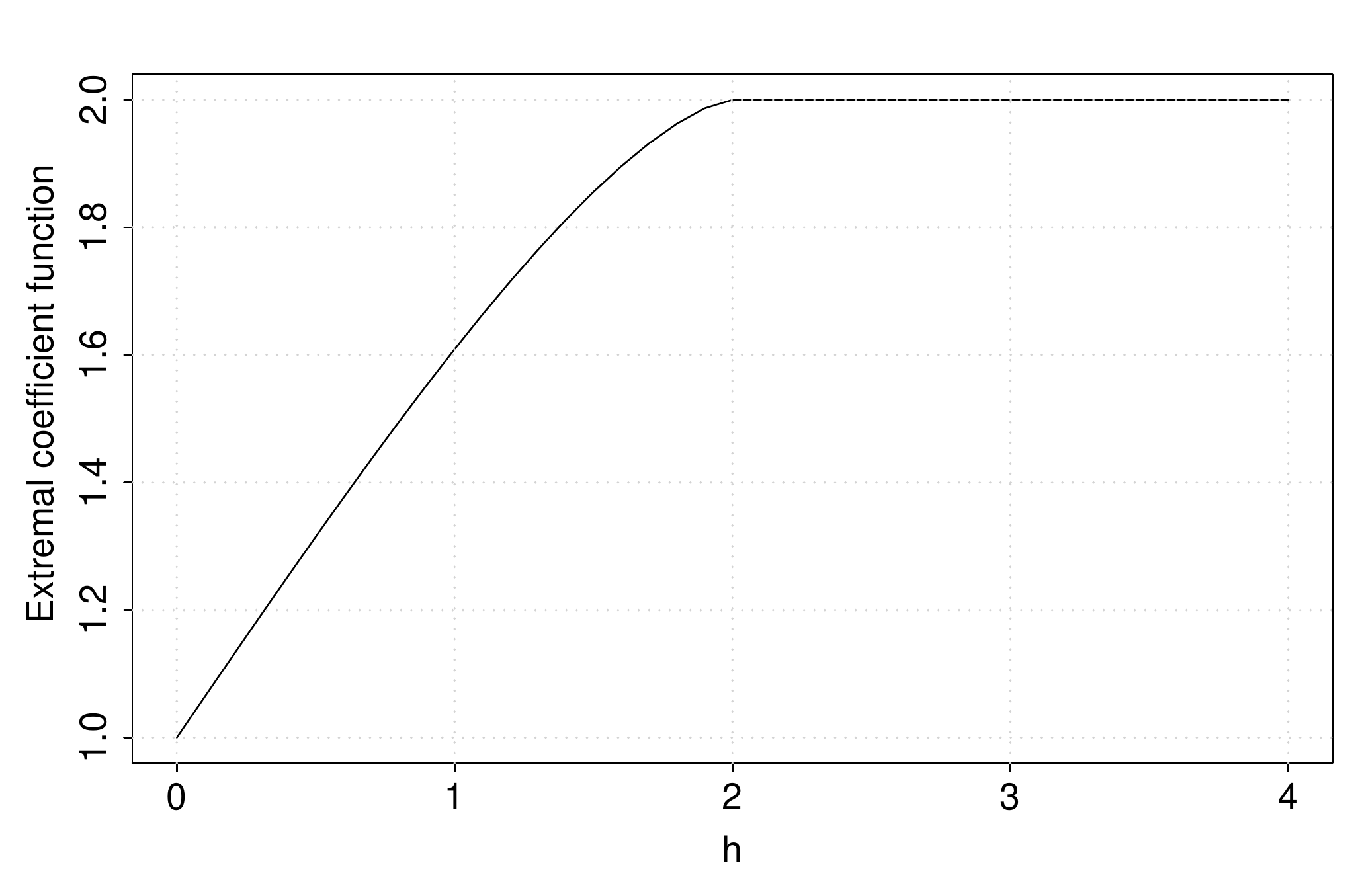}
\center
\caption{Extremal coefficient function of the tube process.}
\label{Extremal coefficient max-stable process with tubes}
\end{figure}

\begin{Def}[Tube process] The tube process is a M2 process where the deterministic shape function is written $f(\Mb{y})=h_b \ \mathds{I}_{ \{ \| \Mb{y} \| < R_b \} }, \Mb{y} \in \mathds{R}^2$, with $R_b>0$ and $h_b=\frac{1}{\pi R_b^2}$.
\end{Def}
The density $f$ is a tube with height $h_b$ and radius $R_b$ centered at point $\Mb{0}$. This process is simple and its bivariate distribution function is given in the next proposition.
\begin{Prop}
\label{Prop_Bivariate_Distr_Function_Tube}
Let $\Mb{x}_1, \Mb{x}_2 \in \mathds{R}^2$, $h=\| \Mb{x}_1-\Mb{x}_2 \|$ and $z_1, z_2>0$.
The bivariate distribution function of the tube process is given by
$$
- \log ( \mathds{P}(Z(\Mb{x}_1) \leq z_1, Z(\Mb{x}_2) \leq z_2) )
= \left \{ 
\begin{array}{ll}
\frac{h_b}{z_1}[\pi R_b^2 -A_{int}(h)]+\frac{1}{z_2} & \mbox{ if } z_2 \leq z_1, \\
\frac{1}{z_1} + \frac{h_b}{z_2}[\pi R_b^2 -A_{int}(h)] & \mbox{ if } z_2 > z_1, \\
\end{array}
\right.
$$
where 
$$
A_{int}(h)=
\left\{
\begin{array}{l l}
2 \bigg(R_b^2 \arcsin \left( \frac{ \sqrt{4R_b^2-h^2}}{ 2 R_b}  \right) - \frac{h}{4} \sqrt{4 R_b^2-h^2} \bigg) & \mbox{ if } \  h \leq 2 R_b,  \\
0 & \mbox{ if } \  h > 2 R_b.
\end{array}
\right.
$$
\end{Prop}
The extremal coefficient function, which is isotropic, is given in the next corollary and depicted in Figure \ref{Extremal coefficient max-stable process with tubes} for $R_b=1$.
\begin{Corr}
The extremal coefficient function of the tube process is given by
\label{Maxstable_Tubes}
$$ \Theta(h)=
\left\{
\begin{array}{l l}
2 \left[1- h_b \bigg(R_b^2 \arcsin \left( \dfrac{ \sqrt{4R_b^2-h^2}}{ 2 R_b}  \right) - \dfrac{h}{4} \sqrt{4R_b^2-h^2} \bigg)\right] & \mbox{ if } \ 0 \leq h \leq 2 R_b,  \\
2 & \mbox{ if } \  h > 2 R_b.
\end{array}
\right.$$
\end{Corr}
\begin{figure}[h!]
\center
\includegraphics[scale=0.52]{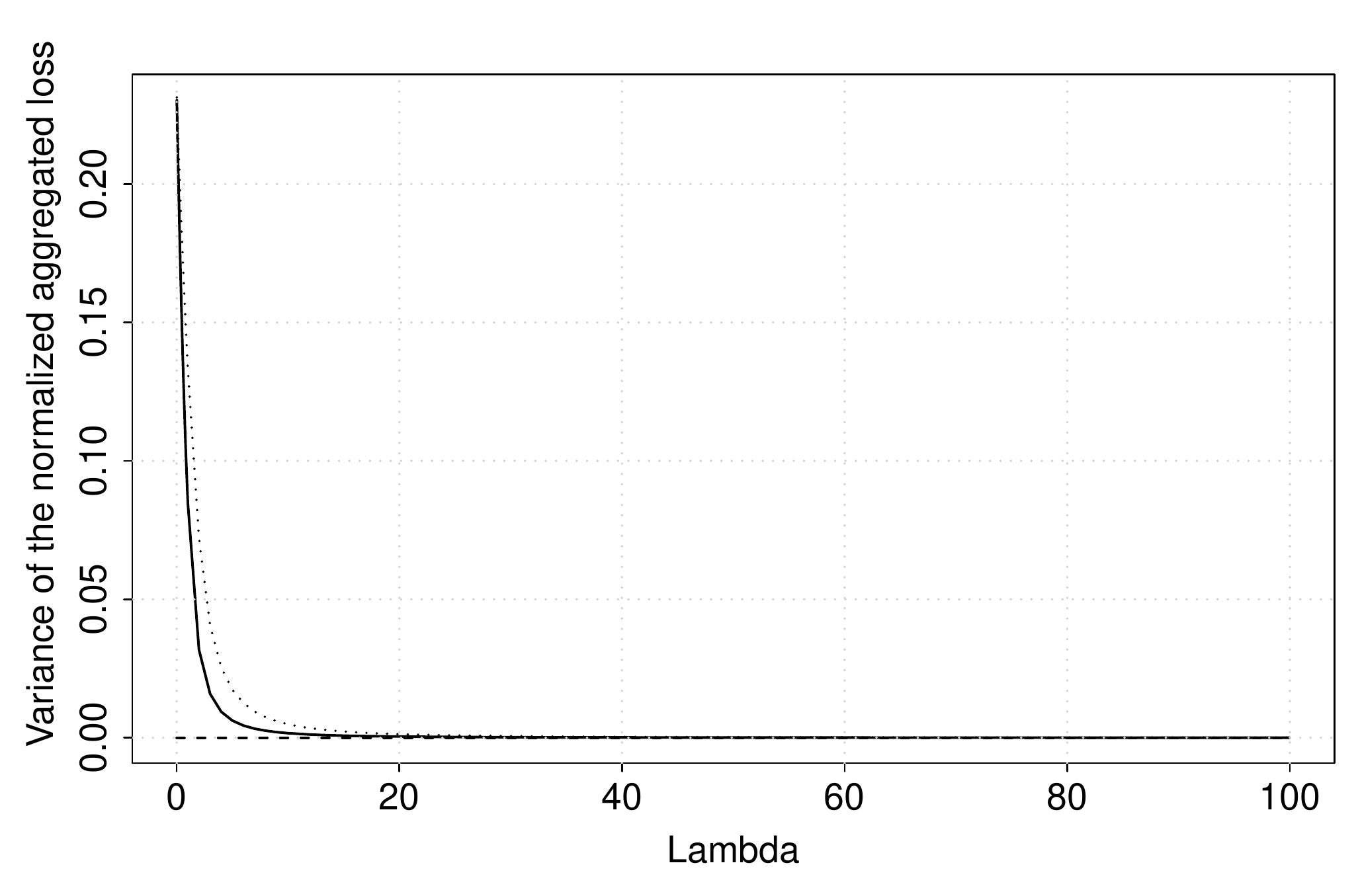}
\caption{The solid (respectively dotted) line depicts the evolution of $\mathcal{R}_2(\lambda A)$ with respect to $\lambda$ in the case of the tube process, where $A$ is a disk (respectively a square). The dashed line represents the limiting risk measure.}
\label{R2 for max-stable process tubes}
\end{figure}
We can easily show that for all $h>0$, the function $\lambda \mapsto \theta(\lambda h)$ is increasing for $\lambda \leq \dfrac{2 R_b}{h}$ and constant above that. Thus, \eqref{Risk_Measure_Disk} and \eqref{Risk_Measure_Square} yield that the function $\lambda \mapsto \mathcal{R}_2(\lambda A)$ is decreasing for $A$ being a disk or a square. Hence, there is spatial diversification. An interesting property stems from the fact that $\Theta(h)$ reaches 2 whenever $h \geq 2 R_b$, meaning that there is spatial independence at a finite distance. This explains intuitively why the spatial diversification is faster than in the case of the previously introduced models, as can be observed in Figure \ref{R2 for max-stable process tubes} (obtained with $R_b=1$).
Obviously, $\lim_{\lambda \to \infty} \Theta(\lambda)=2$, and Corollary \ref{Chapriskmeasures_Corr_Asymptotic_Independence} yields  $\lim_{\lambda \to \infty} \mathcal{R}_2(\lambda A)=0$. Moreover, due to the spatial independence at finite distance, we have the following proposition.
\begin{Prop}
\label{Proposition_Variance_Tubes}
Let $\{ Z(\Mb{x}) \}_{\Mb{x} \in \mathds{R}^2}$ be the tube process and $A \in \mathcal{A}$. We have, for all $\lambda \geq 1$, $\lim_{R_b \to 0} \mathcal{R}_2(\lambda A) = \mathcal{R}_2(A) =0.$
\end{Prop}
\noindent The limit process arising as $R_b$ tends to $0$ corresponds to the case of asymptotic independence whatever the distance.

\medskip

\begin{figure}[h!]
\center
\includegraphics[scale=0.56]{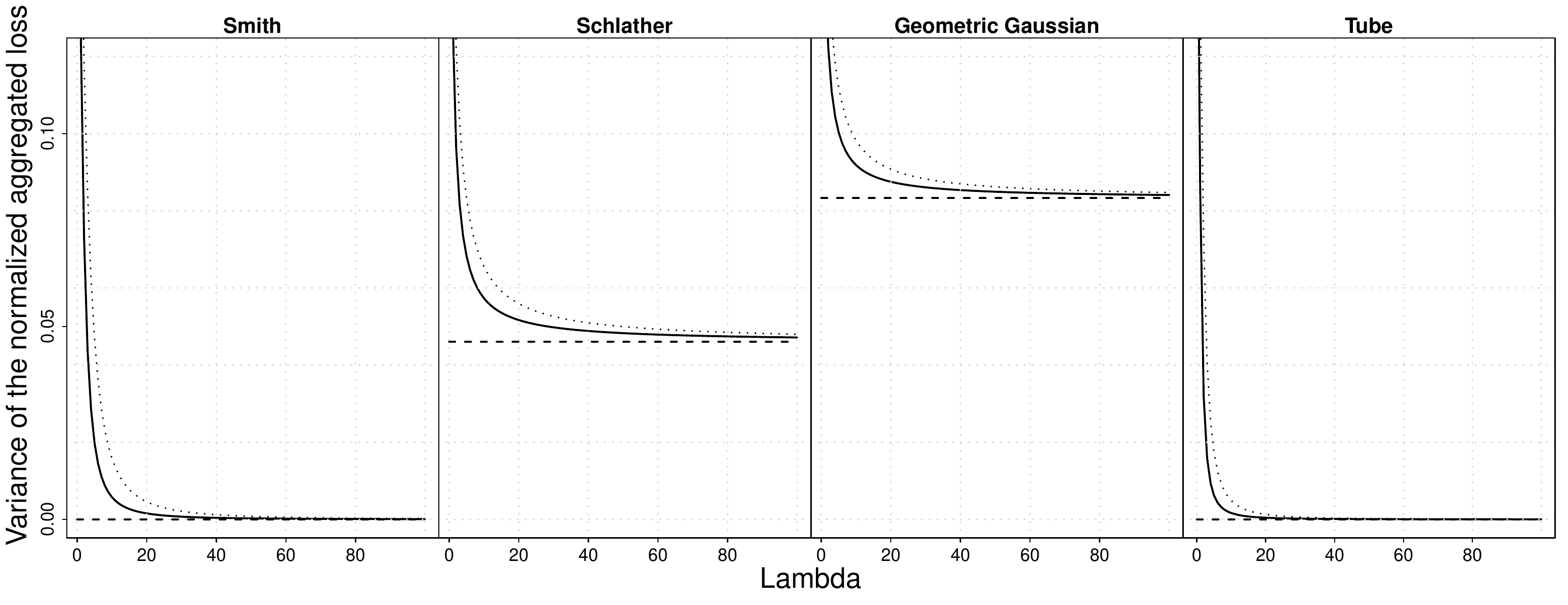}
\caption{Each panel corresponds to a different max-stable model. The solid (respectively dotted) line depicts the evolution of $\mathcal{R}_2(\lambda A)$ with respect to $\lambda$, where $A$ is a disk (respectively a square). The dashed line represents the limiting risk measure.}
\label{R2 Comparison}
\end{figure}

By way of summary, a comparison of the four models considered (Smith, Schlather, geometric Gaussian and tube) is provided in Figure \ref{R2 Comparison}. In the cases of the Schlather and the geometric Gaussian models, a Cauchy correlation function has been used.
These different processes show a large variety of behavior, both in terms of speed and ``completeness'' of spatial diversification. The Brown-Resnick model itself includes various types of behavior. The two particular cases considered here (the Smith model and the geometric Gaussian model) are very different. In all cases, the behavior of the function $\lambda \mapsto \mathcal{R}_2(\lambda A)$ is very similar in the case of the disk and the square. The only difference consists in the fact that the spatial diversification is slightly slower for the square. For an insurance company, the characteristic dimension of the geographical area required to reach a specified low variance level depends to a large extent on the type of max-stable process driving extreme events and on the geometric shape of $A$.

\begin{Rq}
It is important to note that our analysis has been carried out with a standardized characteristic distance of spatial dependence: the eigenvalues of $\Sigma$, the range parameter $c_1$ and the radius of the tube $R_b$ are equal to 1. However, in a specific case study, the characteristic dimension of the area required to reach a given level of variance depends on the real values of $\Sigma$, $c_1$ or $R_b$ on the region of interest.
\end{Rq}

\subsubsection{Axioms} 

The following theorem shows that the spatial risk measure $\mathcal{R}_2$ satisfies, at least partially, the different axioms introduced in Definition \ref{Chapriskmeasures_Def_Axiomatic} under some conditions.
\begin{Th}
\label{Variance_Properties}
1. For any stationary max-stable process $\{ Z(\Mb{x}) \}_{\Mb{x} \in \mathds{R}^2}$, the spatial risk measure $\mathcal{R}_2$ satisfies the axiom of spatial invariance under translation. In particular, this is true for the Smith, Schlather, Brown-Resnick, geometric Gaussian and tube processes.  \\
\noindent 2. For any simple and stationary max-stable process $\{ Z(\Mb{x}) \}_{\Mb{x} \in \mathds{R}^2}$ having an isotropic and non-decreasing extremal coefficient function, $\mathcal{R}_2$ satisfies the axiom of spatial sub-additivity when the two regions are both a disk or a square. In particular, this is true for the Smith process having a covariance matrix proportional to the identity, the Schlather and geometric Gaussian processes having an isotropic and non-increasing correlation function, the Brown-Resnick process having an isotropic and non-decreasing semivariogram and the tube process. \\
\noindent 3. Let $\{ Z(\Mb{x}) \}_{\Mb{x} \in \mathds{R}^2}$ be a simple max-stable process having an extremal coefficient function $\Theta$ such that for all $\Mb{x}, \Mb{y} \in \mathds{R}^2$, $\Theta(\Mb{x}, \Mb{y})$ only depends on $\Mb{x}-\Mb{y}$. Thus, $\Theta(\Mb{x}, \Mb{y})$ is written $\Theta(\Mb{x}-\Mb{y})$.
If 
\Beq
\label{Eq_Condition_Thm_3}
0 < \int_{\mathds{R}^2} (2-\Theta(\Mb{x})) \ d\Mb{x} < \infty,
\Eeq
then the spatial risk measure $\mathcal{R}_2$ satisfies the axiom of asymptotic spatial homogeneity of order $-2$, with 
$$ K_1=0 \quad \mbox{and} \quad K_2=\dfrac{\sigma^2}{|A|},$$
where, for $u>0$,
\Beq
\label{Eq_Def_Sigma}
\sigma=\sqrt{ \int_{\mathds{R}^2} \left[ \exp \left( -\frac{\Theta(\Mb{x})}{u} \right) - \exp \left( -\frac{2}{u} \right) \right] d\Mb{x}}.
\Eeq
In particular, this is true for the Smith process, the Brown-Resnick process with a semivariogram satisfying $\gamma( \Mb{h} ) = \eta \| \Mb{h} \|^a$, $\Mb{h} \in \mathds{R}^2$, where $\eta>0$ and $a \in (0, 2]$, and the tube process.
\end{Th}

\begin{Rq}
\label{Rq_Generalization_1}
It turns out that the axiom of asymptotic spatial homogeneity of order $-2$ of spatial risk measures associated with the variance is satisfied for more general cost processes than the one considered in this section. Indeed, let Cov denote the covariance and $\{ C_P(\Mb{x}) \}_{\Mb{x} \in \mathds{R}^2}$ be a cost process having locally integrable sample paths and such that, for all $\Mb{x}, \Mb{y} \in \mathds{R}^2$, $\mbox{Cov}(C_P(\Mb{x}), C_P(\Mb{y}))$ only depends on $\Mb{x}-\Mb{y}$, and $\sigma_g^2=\displaystyle \int_{\mathds{R}^2} \mbox{Cov}(C_P(\Mb{0}), C_P(\Mb{x})) \ d\Mb{x} \in (0, \infty)$.
From the proof of Theorem \ref{Variance_Properties}, Bullet 3, we easily see that 
the quantity $\mbox{Var} \left( \dfrac{1}{|A|} \displaystyle \int_A C_P(\Mb{x}) \ d\Mb{x} \right)$ satisfies the axiom of asymptotic spatial homogeneity of order $-2$ with $K_1=0$ and $K_2=\dfrac{\sigma_g^2}{|A|}$.
\end{Rq}

\subsection{Central limit theorem}

Having shown that, under some conditions, $\mathcal{R}_2$ satisfies the axiom of asymptotic spatial homogeneity of order $-2$, the question whether $L_N(\lambda A)$ (for $A$ being a convex in $\mathcal{A}$) satisfies a central limit theorem arises pretty naturally. Using a result by \cite{spodarev2014limit}, we can derive the following result.
\begin{Th} 
\label{TCL_Temperature}
Let $\{ Z(\Mb{x}) \}_{\Mb{x} \in \mathds{R}^2}$ be a simple, stationary and measurable max-stable process having $\Theta$ as extremal coefficient function. If 
\Beq
\label{Eq_Condition_Thm_4}
\int_{\mathds{R}^2} (2-\Theta(\Mb{x})) \ d\Mb{x} < \infty,
\Eeq
we have, for all convex $A \in \mathcal{A}$, that
\Beq
\label{Eq_Convergence_TCL}
\lambda \left( L_N(\lambda A) - m \right) \overset{d}{\to} \mathcal{N} \left( 0, \frac{\sigma^2}{|A|} \right), \mbox{ for } \lambda \to \infty,
\Eeq
where
$$m=1 - \exp \left( -\dfrac{1}{u} \right).$$
In particular, this is true for the Smith process, the Brown-Resnick process with a semivariogram satisfying $\gamma( \Mb{h} ) = \eta \| \Mb{h} \|^a$, $\Mb{h} \in \mathds{R}^2$, where $\eta>0$ and $a \in (0, 2]$, and the tube process.
\end{Th}
It should be noted that there may be mixing max-stable processes for which \eqref{Eq_Condition_Thm_4} is not satisfied. This is the case if $\Theta(\Mb{x})$ does not converge fast enough to $2$ as $\| \Mb{x} \|$ tends to infinity.

\subsection{The Value-at-Risk}
\label{Chapriskmeasures_Subsec_VaR}

Here we focus on $\mathcal{R}_{3, \alpha}(A)=\mbox{VaR}_{\alpha} ( L_N(A) )$, $\alpha \in (0,1)$, where, for a random variable $X$ with distribution function $F$, $\mbox{VaR}_{\alpha}(X)=\inf \{ x \in \mathds{R}: F(x) \geq \alpha \}$. The risk measure VaR is of interest for insurance companies, and this mainly from a regulatory point of view \cite[see e.g.][]{QRM2015}. We adopt an approach similar to the case of the variance, giving first some insights about the function $\lambda \mapsto \mathcal{R}_{3, \alpha}(\lambda A)$ before discussing whether $\mathcal{R}_{3, \alpha}$ satisfies the axioms of Definition \ref{Chapriskmeasures_Def_Axiomatic}.

\subsubsection{Study of $\mathcal{R}_{3, \alpha}(\lambda A), \lambda >0$}

Since deriving a tractable formula for the VaR of $L_N(\lambda A)$ is very difficult, in the case of max-stable processes which can be simulated, we propose to evaluate it as follows. We simulate the process under consideration on a grid containing the locations $\Mb{x}_m \in \lambda A, m=1, \dots, \lambda^2 M$, where $M$ is the number of sites in region $A$. Then, one realization of $L_N(\lambda A)$ can be approximated by a Riemann sum. By generating a number $S$ of independent approximated replications of the random variable $L_N(\lambda A)$, an approximation of its distribution function is obtained. Finally, taking the empirical quantile of this distribution at level $\alpha$ yields an approximation of $\mathcal{R}_{3, \alpha}$.

\begin{figure}[h!]
\center
\includegraphics[scale=0.52]{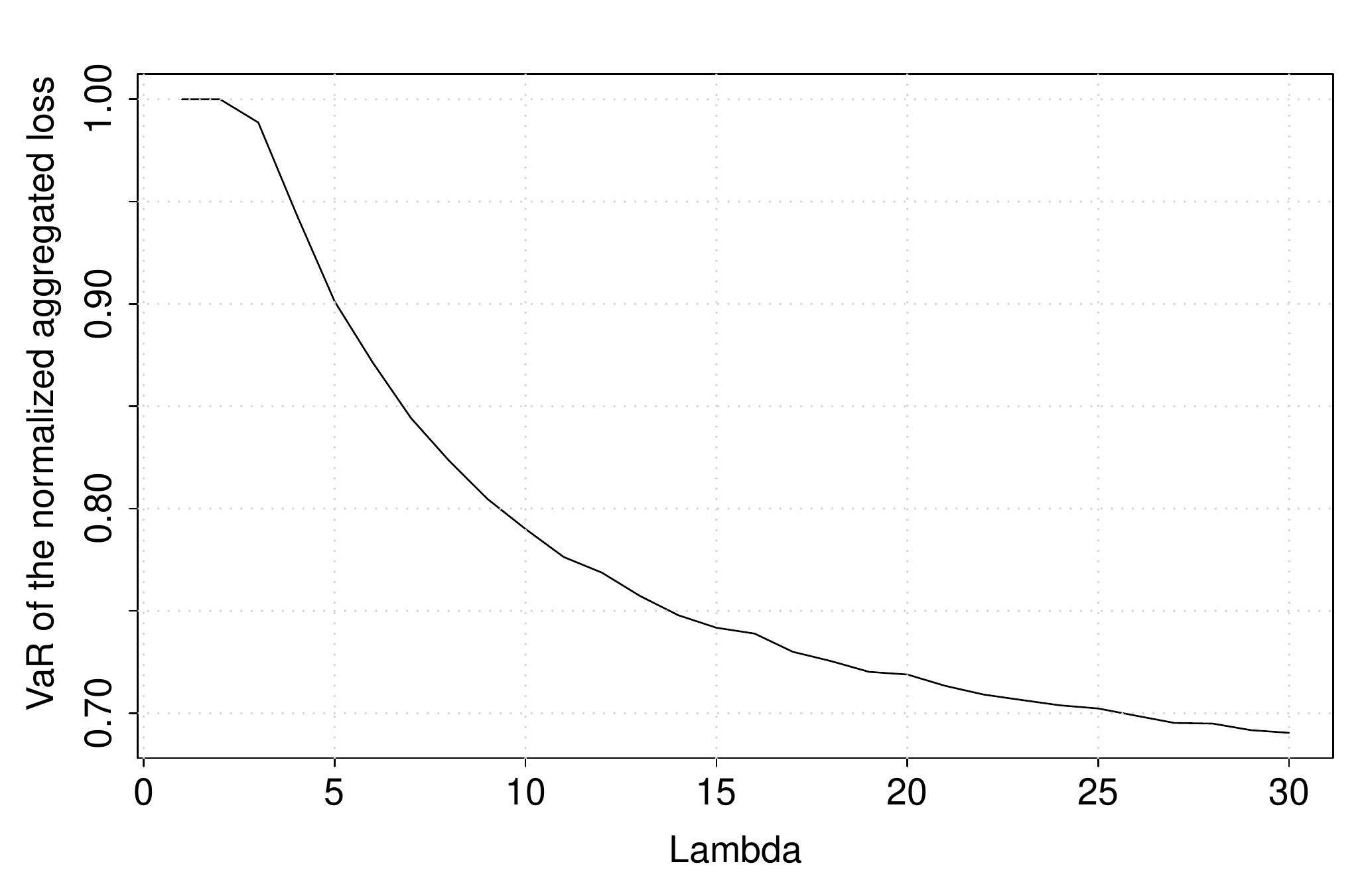}
\caption{$\mathcal{R}_{3,0.9}(\lambda A)$ with respect to $\lambda$ in the case of the Smith process, where $A$ is a square of side 1.}
\label{R3 for Smith process}
\end{figure}

To illustrate this procedure, the evolution of $\mathcal{R}_{3,0.9}(\lambda A)$ with respect to $\lambda$ in the case of the Smith model is depicted in Figure \ref{R3 for Smith process}, for $M=49$ and $S=10000$. The sites are located on a regular grid.
Of course, the larger $M$ and $S$, the smoother the curve. Figure 7 seems to indicate spatial diversification.

\subsubsection{Axioms}

Although we do not have an explicit formula for $\mathcal{R}_{3, \alpha}(\lambda A)$, by taking advantage of Theorem \ref{TCL_Temperature}, we know the asymptotic behavior (when $\lambda \to \infty$) of $\mathcal{R}_{3, \alpha}(\lambda A)$ under some conditions. This result is part of the following theorem.

\begin{Th}
\label{Th_Asymptotic_Homogeneity_VaR}
1. For any stationary max-stable process $\{ Z(\Mb{x}) \}_{\Mb{x} \in \mathds{R}^2}$, the spatial risk measure $\mathcal{R}_{3, \alpha}$ satisfies the axiom of spatial invariance under translation. In particular, this is true for the Smith, Schlather, Brown-Resnick, geometric Gaussian and tube processes.  \\
\noindent 2. For any simple, stationary and measurable max-stable proces $\{ Z(\Mb{x}) \}_{\Mb{x} \in \mathds{R}^2}$ with extremal coefficient function $\Theta$ such that \eqref{Eq_Condition_Thm_4} is satisfied, the spatial risk measure $\mathcal{R}_{3, \alpha}$ satisfies, for $\alpha \neq 0.5$, the axiom of asymptotic spatial homogeneity of order -1, with 
$$ K_1=1 - \exp \left( -\dfrac{1}{u} \right) \quad \mbox{and} \quad K_2= \dfrac{\sigma \ q_{\alpha}}{\sqrt{|A|}},$$
where $q_{\alpha}$ is the quantile at the level $\alpha$ of the standard Gaussian distribution. In particular, this is true for the Smith process, the Brown-Resnick process with a semivariogram satisfying $\gamma( \Mb{h} ) = \eta \| \Mb{h} \|^a$, $\Mb{h} \in \mathds{R}^2$, where $\eta>0$ and $a \in (0,2]$, and the tube process.
\end{Th}
Since we do not have an explicit formula for $\mathcal{R}_{3, \alpha}(\lambda A)$ for $\lambda < \infty$, there is nothing we can say about spatial sub-additivity. However, the previous numerical study suggests that, in the case of the Smith model, the function $\lambda \mapsto \mathcal{R}_{3, \alpha}(\lambda A)$ is non-increasing for $A$ being a square, which would imply spatial sub-additivity when the two regions are both a square.

\begin{Rq}
\label{Rq_Generalization_2}
From the proof of Theorem \ref{Th_Asymptotic_Homogeneity_VaR}, Bullet 2, we see that asymptotic spatial homogeneity of order $-1$ of spatial risk measures associated with VaR is satisfied as soon as there is a central limit theorem for the normalized spatially aggregated loss. The existence of such a theorem can be linked to the mixing properties of the cost process \citep[see e.g.][Theorem 1]{gorodetskii1987moment}.
\end{Rq}

\subsection{Practical considerations}

In our analysis, we have most often considered max-stable processes with standard Fr\'echet margins although in practice, the univariate marginal distributions of max-stable processes are Generalized Extreme Value distributions, involving three parameters: $\mu  \in \mathds{R}$, $\sigma>0$ and $\xi \in \mathds{R}$. These distributions are denoted $\mbox{GEV}(\mu, \sigma, \xi)$. For a random variable $Y$, we have that
$Y \sim \mbox{GEV}(\mu, \sigma, \xi) \mbox{ iff } Z=\left[ 1 + \xi \left( \frac{Y-\mu}{\sigma} \right) \right]^{\frac{1}{\xi}} \sim \mbox{standard Fr\'echet}$. 
Hence, it is easy to show that, for $u_1 \in \mathds{R}$,
$Y>u_1 \mbox{ iff } Z> u$, where  
\Beq
\label{Eq_Transformation_u1_u}
u = 
\left \{
\begin{array}{ll}
\left[ 1+ \xi \left( \frac{u_1-\mu}{\sigma} \right) \right]^{\frac{1}{\xi}} & \mbox{ if } \xi \neq 0, \\
\exp \left( \frac{u_1-\mu}{\sigma} \right) & \mbox{ if } \xi=0. 
\end{array}
\right.
\Eeq
The parameters $\mu$, $\sigma$ and $\xi$ are such that $u>0$. This means that the general case of identical GEV margins reduces to the standard Fr\'echet case we have studied. A concrete application of the tools proposed in this section involves, in the case of identical margins, the following steps:
\Ben
\item Specify the threshold $u_1$ of interest, for instance $u_1=30 ^{\circ} C$.
\item Jointly fit to historical data the marginal parameters $\mu$, $\sigma$ and $\xi$ as well as the dependence parameters for different models (e.g. Smith, Schlather, Brown-Resnick, geometric Gaussian, tube), using, for instance, the composite likelihood approach  \citep[see e.g.][]{padoan2010likelihood}.
\item Carry out model selection based, for instance, on the Akaike Information Criterion \citep{akaike1974new} or the likelihood ratio statistic \citep{davison2003statistical}.
\item Transform $u_1$ in $u$ by plugging the estimates of $\mu$, $\sigma$ and $\xi$ in \eqref{Eq_Transformation_u1_u}. 
\item Study different spatial risk measures as explained, considering the threshold $u$ computed in the previous step.
\Een
In some situations, the marginal parameters may vary across space, meaning, using \eqref{Eq_Transformation_u1_u}, that the threshold also varies across space and is now of the form $u(\Mb{x})$, $\Mb{x} \in \mathds{R}^2$. In this case, $\mathcal{R}_2(\lambda A)$ involves the full expression of the bivariate distribution function and not only the extremal coefficient, leading to more complex (but still numerically computable) expressions than \eqref{Chapriskmeasures_Eq_Stationary_Risk_Measure_Generalcase}, \eqref{Risk_Measure_Disk} and \eqref{Risk_Measure_Square}. If the dependence parameters also vary across space, the complexity becomes even higher.

\section{Conclusion}
\label{Chap_RislMeasures_Sec_Conclusion}

The first part of the paper introduces a new notion of spatial risk measure, based on the normalized spatially aggregated loss, and proposes a set of axioms adapted to the spatial context. These axioms appear natural for stationary cost processes. Contrary to the classical literature relating to risk measures, our axiomatic approach aims at quantifying the sensitivity of the risk measurement with respect to the space variable. 
Our motivation comes from insurance applications where losses are due to hazards having a spatial extent, typically such as environmental hazards. In order to develop concrete examples of spatial risk measures, we propose a general model that maps the process of the environmental variable into a cost, via a damage function. 

The second part of the paper considers a specific model for the cost, where the environmental process is max-stable and the damage function is the indicator of threshold exceedances. Some spatial risk measures based on this model are studied from a theoretical viewpoint. Under some conditions mainly on the underlying max-stable process, they satisfy, at least partially, the axioms proposed. This feature is true for more general cost processes, as shown especially in Remarks \ref{Rq_Generalization_1} and \ref{Rq_Generalization_2}.
Characterizing all univariate risk measures and processes allowing the proposed axioms to be satisfied would be relevant and useful.

In the case of non-stationary cost processes, more adapted axioms should be proposed, depending of course on the type of non-stationarity. If the cost process can be assumed to be piecewise stationary, the concepts introduced here can be independently applied to each sub-region on which stationarity may be assumed. Asymptotic spatial homogeneity naturally requires sufficiently large sub-regions.

\section*{Acknowledgements} 

The author wishes to thank Paul Embrechts and Christian-Yann Robert for their many useful comments and suggestions. He also acknowledges Romain Biard, Cl\'ement Dombry, Nabil Kazi-Tani, Pierre Ribereau, Gennady Samorodnitsky, an anonymous Associate Editor, two anonymous referees, and the participants in the financial econometrics seminar at CREST, the CFE-ERCIM 2013 conference, the PEPER Workshop 2013, the RiskLab seminar at ETH Zurich, the Lyon-Lausanne seminar, the seminar at the ``Laboratoire de Math\'ematiques de Besan\c{c}on'', the CNRM-GAME seminar at ``M\'et\'eo France'', the Uncertainty, Risk and Robustness Workshop at Columbia University, the Mathematical Foundations of Heavy Tailed Analysis Workshop at the University of Copenhagen, the Oberseminar Finanz- und Versicherungsmathematik at the Technical University of Munich and the Workshop on ``Understanding the Risks of Extreme Events: Challenges in Risk Management of Natural Disasters'' at ETH Zurich. The author would like to thank the MIRACCLE-GICC project, RiskLab at ETH Zurich and the Swiss Finance Institute for financial support.

\newpage
\appendix
\section{Appendix: Proofs}

\subsection{For Theorem \ref{Prop_Stationary_Risk_Measure_Generalcase}}
\begin{proof}

We have, for all $A \in \mathcal{A}$, that
\begin{align*}
\mathds{E} \left( [ L(A) ]^2 \right) = \mathds{E}\left( \left( \int_A \mathds{I}_{\{Z(\Mb{x})>u\}} \ d\Mb{x}  \right)^2 \right) &= \mathds{E} \left( \int_A \mathds{I}_{\{Z(\Mb{x})>u\}} \ d\Mb{x} \ \int_A \mathds{I}_{\{Z(\Mb{y})>u\}} \ d\Mb{y}  \right)
\\& = \mathds{E} \left( \int_A \int_A \mathds{I}_{\{Z(\Mb{x})>u, Z(\Mb{y})>u\}} \ d\Mb{x}\  d\Mb{y} \right)
\\& = \int_A \int_A \mathds{P}(Z(\Mb{x})>u, Z(\Mb{y})>u)\  d\Mb{x}\  d\Mb{y}.
\end{align*}
Moreover, using the fact that $Z$ has standard Fr\'echet margins, we have
$$ \mathds{P} ( Z(\Mb{x})>u, Z(\Mb{y})>u)= 1+ \mathds{P}(Z(\Mb{x}) \leq u, Z(\Mb{y}) \leq u) -2 \exp \left( -\frac{1}{u} \right), \mbox{ yielding } $$
$$ \mathds{E} \left( (L(A))^2 \right) = |A|^2 \left[ 1-2 \exp \left( -\frac{1}{u} \right) \right] + \int_A \int_A \exp\left( -\frac{\Theta(\Mb{x}, \Mb{y})}{u} \right) \  d\Mb{x}\  d\Mb{y}.$$
Hence,
\begin{align*}
\mathcal{R}_2(A) &= \frac{1}{|A|^2} \Bigg( |A|^2 \left[ 1-2 \exp \left( -\frac{1}{u} \right) \right] + \int_A \int_A \exp\left( -\frac{\Theta(\Mb{x}, \Mb{y})}{u} \right) \  d\Mb{x}\  d\Mb{y} 
\\& \quad \ -  |A|^2  \left[ 1 - \exp \left( -\frac{1}{u} \right) \right]^2 \Bigg)
\\& =  \frac{1}{|A|^2} \int_A  \int_A \left[ \exp\left( -\frac{\Theta(\Mb{x}, \Mb{y})}{u} \right) - \exp \left( -\frac{2}{u} \right) \right] \  d\Mb{x}\  d\Mb{y}.
\end{align*}
The result is obtained by replacing $A$ with $\lambda A$.
\end{proof}

\subsection{For Corollary \ref{Prop_Stationary_Risk_Measure}} 

\begin{proof}
1. Let $A$ be a disk with radius $R$. For any integrable function $g: (\Mb{x}, \Mb{y}) \in \mathds{R}^2 \times \mathds{R}^2 \mapsto \mathds{R}$ depending only on $ h= \| \Mb{x}-\Mb{y} \|$, it is easy to show that, for all $\lambda >0$,
\Beq
\label{Chapriskmeasures_Transformation_Distance}
\int_{\lambda A} \int_{\lambda A} g(\Mb{x},\Mb{y}) \ d\Mb{x}\  d\Mb{y} = \lambda^4 |A|^2 \int_{0}^{2 \lambda R} \ f_d(h,\lambda R)\  g(h)\  dh.
\Eeq
Moreover, we have \citep[see e.g.][Section 4.3.1]{moltchanov2012distance} that
$$
f_d(h,R)=\frac{2h}{R^2} \left( \frac{2}{\pi} \arccos \left( \frac{h}{2 R} \right) - \frac{h}{\pi R} \sqrt{1-\frac{h^2}{4 R^2}}   \right), \mbox{ for } 0 \leq h \leq 2R,
$$
yielding, for all $\lambda>0$, $f_d(\lambda h,\lambda R)= \dfrac{1}{\lambda} \ f_d(h,R).$
Therefore, by the change of variable $h_d=\dfrac{h}{\lambda}$, we obtain
\Beq
\label{Eq_Variable_Change}
\int_{0}^{2 \lambda R} f_d(h,\lambda R) \ g(h)\  dh
= \int_{0}^{2 R} f_d(h_d, R ) \ g(\lambda h_d) \  dh_d.
\Eeq
Applying \eqref{Chapriskmeasures_Transformation_Distance} and \eqref{Eq_Variable_Change} to $g(\Mb{x},\Mb{y})=\exp \left( -\dfrac{\Theta(\Mb{x}, \Mb{y})}{u} \right)=\exp \left( -\dfrac{\Theta( \| \Mb{x}-\Mb{y} \| )}{u} \right)$ (since $\Theta$ is assumed to be isotropic), and using \eqref{Chapriskmeasures_Eq_Stationary_Risk_Measure_Generalcase}, we obtain
$$
\mathcal{R}_2 (\lambda A) = - \exp \left( -\frac{2}{u} \right) + \int_{0}^{2 R} f_d(h,R) \ \exp \left( -\frac{\Theta(\lambda h)}{u} \right) \  dh.
$$

\medskip

\noindent 2. The same reasoning yields the result for $A$ being a square. The distribution function of the Euclidean distance between two points independently
and uniformly distributed on a square with side $R$ is written \citep[see e.g.][Section 4.3.2]{moltchanov2012distance}
\begin{align*}
& \quad \ F_s(h,R) 
\\&=
\left \{
\begin{array}{ll}
\dfrac{\pi h^2}{R^2}-\dfrac{8 h^3}{3 R^3}+\dfrac{h^4}{2R^4} & \mbox{ if } h \in [0,R], \\
\dfrac{1}{3}-2b-\dfrac{b^2}{2}+\dfrac{2}{3} \sqrt{(b-1)^3}+ 2 \sqrt{b-1} + 2b \sqrt{b-1}+2b \arcsin \left( \dfrac{2-b}{b} \right) & \mbox{ if } h \in  [R, R \sqrt{2}],
\end{array}
\right.
\end{align*}
where $b=\dfrac{h^2}{R^2}$. Hence, for $h \in [0,R]$, the density is 
$$ f_s(h,R)=\frac{2 \pi h}{R^2} - \frac{8 h^2}{R^3} + \frac{2 h^3}{R^4}.$$
For $h \in  [R, R \sqrt{2}]$, we obtain
$$
f_s(h,R)=\left( -2-b+ 3 \sqrt{b-1} + \frac{b+1}{\sqrt{b-1}} +2 \arcsin \left( \frac{2-b}{b} \right) - \frac{4}{b \sqrt{1-\frac{(2-b)^2}{b^2}}}\right) \frac{2h}{R^2}.
$$

\medskip

\noindent 3. Since $f_d$ is continuous on the compact set $[0, 2R]$, it is bounded. Thus, using Lebesgue's dominated convergence Theorem and the fact that $f_d$ is a density, we obtain
\begin{align*}
\lim_{\lambda \to \infty} \int_{0}^{2R} f_d(h,R) \  \exp \left( -\frac{\Theta(\lambda h)}{u} \right) \  dh &= \int_{0}^{2R} f_d(h,R) \  \exp \left( -\frac{ \lim_{\lambda \to \infty} \Theta(\lambda)}{u} \right) \  dh 
\\& = \exp \left( -\frac{ \lim_{\lambda \to \infty} \Theta(\lambda)}{u} \right),
\end{align*}
which gives the result when $A$ is a disk. The argument is the same in the case of the square.
\end{proof}

\subsection{For Proposition \ref{Prop_Bivariate_Distr_Function_Tube}}

\begin{proof}
Since the $(U_i, \Mb{S}_i)_i$ are the points of a Poisson point process on $(0,\infty) \times \mathds{R}^2$ with intensity $u^{-2} du \times d\Mb{s}$, we have, for all $z_1, z_2 >0$, that
\begin{align}
&- \log ( \mathds{P}(Z(\Mb{x}_1) \leq z_1, Z(\Mb{x}_2) \leq z_2) ) \nonumber
\\& = \int_{\mathds{R}^2} \int_{\min \left \{ \dfrac{z_1}{f(\Mb{x}_1-\Mb{s})},\dfrac{z_2}{f(\Mb{x}_2-\Mb{s})} \right \} }^{\infty} u^{-2} du \ d\Mb{s} \nonumber
\\& = \int_{\mathds{R}^2} \max \left \{ \frac{f(\Mb{x}_1-\Mb{s})}{z_1},\frac{f(\Mb{x}_2-\Mb{s})}{z_2} \right \} \ d\Mb{s} \nonumber
\\& = \int_{\mathds{R}^2} \frac{f(\Mb{x}_1-\Mb{s})}{z_1} \mathds{I}_{ \left \{ \frac{f(\Mb{x}_1-\Mb{s})}{z_1} > \frac{f(\Mb{x}_2-\Mb{s})}{z_2} \right \}} \ d\Mb{s} + \int_{\mathds{R}^2} \frac{f(\Mb{x}_2-\Mb{s})}{z_2} \mathds{I}_{ \left \{ \frac{f(\Mb{x}_2-\Mb{s})}{z_2} \geq \frac{f(\Mb{x}_1-\Mb{s})}{z_1} \right \}} \ d\Mb{s} \nonumber
\\& = \int_{\mathds{R}^2} \frac{f(\Mb{s})}{z_1} \mathds{I}_{\left \{ \frac{f(\Mb{s})}{z_1} > \frac{f(\Mb{s}+\Mb{x}_2-\Mb{x}_1)}{z_2} \right \} } \ d\Mb{s} + \int_{\mathds{R}^2} \frac{f(\Mb{s})}{z_2} \mathds{I}_{ \left \{ \frac{f(\Mb{s})}{z_2} \geq \frac{f(\Mb{s}+\Mb{x}_1-\Mb{x}_2)}{z_1} \right \}} \ d\Mb{s} \nonumber
\\& = \frac{1}{z_1} \mathds{P} \left( \frac{f(\Mb{X})}{z_1} > \frac{f(\Mb{X}+\Mb{x}_2-\Mb{x}_1)}{z_2} \right) 
+ \frac{1}{z_2} \mathds{P} \left( \frac{f(\Mb{X})}{z_2} \geq \frac{f(\Mb{X}+\Mb{x}_1-\Mb{x}_2)}{z_1} \right), 
\label{Chapriskmeasures_Bivariate_Density_Tube}
\end{align}
where $\Mb{X}$ is a random vector having density $f$. \\
\noindent Let $E_1$ be the  event $\left \{ \dfrac{f(\Mb{X})}{z_1} > \dfrac{f(\Mb{X}+\Mb{x}_2-\Mb{x}_1)}{z_2} \right \}= \left \{ z_2 \ \mathds{I}_{ \{ \| \Mb{X} \| \leq R_b \} } > z_1 \ \mathds{I}_{ \{ \| \Mb{X}+\Mb{x}_2-\Mb{x}_1 \| \leq R_b \} } \right \}$. If $z_2 > z_1$,
$E_1 = \left \{  \| \Mb{X} \| \leq R_b \right \}$, giving
$\mathds{P}(E_1)=1.$
Indeed $\Mb{X}$ has density $f$ and then $\| \Mb{X} \| \leq R_b$ almost surely (a.s.).
If $z_1 \geq z_2$,
$E_1 = \left \{  \| \Mb{X} \| \leq R_b \mbox{ and } \| \Mb{X}+\Mb{x}_2-\Mb{x}_1 \| > R_b  \right \} = \left \{ \| \Mb{X}+\Mb{x}_2-\Mb{x}_1 \| > R_b  \right \}
$
since $ \| \Mb{X} \| \leq R_b$ is satisfied a.s.
Therefore,
\begin{align*}
\mathds{P}(E_1)= \int_{\mathds{R}^2}  \mathds{I}_{ \{ \| \Mb{x}+\Mb{x}_2-\Mb{x}_1 \| > R_b \} } \ f(\Mb{x}) \ d\Mb{x} &= h_b \int_{\mathds{R}^2}  \mathds{I}_{ \{ \|  \Mb{x} \| \leq R_b \ \cap \ \| \Mb{x}-(\Mb{x}_1-\Mb{x}_2) \| > R_b \} } \ d\Mb{x} 
\\& = h_b \ [\pi R_b^2 -A_{int}(h)],
\end{align*}
where $A_{int}(h)$ is the area of the intersection between the base of the tube of center $\Mb{0}$ and that of the tube of center $(\Mb{x}_1-\Mb{x}_2)$.
Note that by symmetry, the area of the intersection between the base of the tube of center $\Mb{0}$ and that of the tube of center $(\Mb{x}_2-\Mb{x}_1)$ is also equal to $A_{int}(h)$. \\
\noindent Let $E_2$ be the event
$\left \{ \dfrac{f(\Mb{X})}{z_2} \geq \dfrac{f(\Mb{X}+\Mb{x}_1-\Mb{x}_2)}{z_1} \right \}=\left \{ z_1 \ \mathds{I}_{ \{ \| \Mb{X} \| \leq R_b \} } \geq z_2 \ \mathds{I}_{ \{ \| \Mb{X}+\Mb{x}_1-\Mb{x}_2 \|  \leq R_b \} } \right \}$. If $z_2 > z_1$, $E_2 = \left \{  \| \Mb{X}+\Mb{x}_1-\Mb{x}_2 \| > R_b \right \}$, giving
$ \mathds{P}(E_2)=h_b[\pi R_b^2 -A_{int}(h)]$.
If $z_1 \geq z_2$, $\mathds{P}(E_2)=1$ since $\| \Mb{X} \| \leq R_b$ a.s. \\
Hence, using \eqref{Chapriskmeasures_Bivariate_Density_Tube}, we obtain
\Beq
\label{Bivariate_Distri_Function_Tube}
- \log ( \mathds{P}(Z(\Mb{x}_1) \leq z_1, Z(\Mb{x}_2) \leq z_2) )
= \left \{ 
\begin{array}{ll}
\frac{h_b}{z_1}[\pi R_b^2 -A_{int}(h)]+\frac{1}{z_2} & \mbox{ if } z_2 \leq z_1, \\
\frac{1}{z_1} + \frac{h_b}{z_2}[\pi R_b^2 -A_{int}(h)] & \mbox{ if } z_2 > z_1. \\
\end{array}
\right.
\Eeq

Since we consider the Euclidean norm, the bases of the tubes are circular. Let us then compute the area of the intersection between two discs, respectively with radius $R_b$ and centers $C_1$ and $C_2$ at a distance $h$.
This intersection is not empty iff $h \leq 2 R_b$. We consider this case, represented in the following picture:

\setlength{\unitlength}{1mm} 
\begin{picture}(35,35)(-45,-9)
\put(10,10){\circle{30}}
\put(10,10){\circle*{0.5}}
\put(10,11){$C_1$}
\put(30,10){\circle{30}}
\put(30,10){\circle*{0.5}}
\put(30,11){$C_2$}
\put(19.5,23){$I$}
\end{picture}

\noindent Let us introduce $p= \frac{1}{2} (2R_b+h)$. By using H\'eron's formula, the area of the triangle $IC_1C_2$, denoted by $A_T$, is given by
\Beq
\label{Chapriskmeasures_Heron}
A_T = \sqrt{p(p-R_b)(p-R_b)(p-h)}.
\Eeq
Furthermore, denoting $H$ the height of the triangle $IC_1C_2$, we have
$A_T=\dfrac{h H}{2}$, giving
\Beq
\label{Chapriskmeasures_H}
H=\frac{2 A_T}{h}.
\Eeq
Denote $\alpha$ and $\beta$ the angles $\widehat{IC_1C_2}$ and $\widehat{IC_2C_1}$, respectively. We have
$\sin \alpha=  \sin \beta = \dfrac{H}{R_b}$, yielding, using \eqref{Chapriskmeasures_H},
\Beq
\label{Chapriskmeasures_alpha}
\alpha =\beta = \arcsin \left( \frac{2 A_T}{h R_b} \right).
\Eeq
Let $S$ be the area of the angular sectors delimited respectively by the angles $\alpha$ and $\beta$.
We have
\Beq
\label{Chapriskmeasres_S}
S = \frac{\alpha R_b^2}{2}.
\Eeq
Combining \eqref{Chapriskmeasures_Heron}, \eqref{Chapriskmeasures_alpha} and \eqref{Chapriskmeasres_S}, we obtain, for $h \leq 2 R_b$,
$$
A_{int}(h)= 2 (2S-A_T)
 = 2 \bigg( R_b^2 \arcsin \left( \frac{2 \sqrt{p(p-R_b)^2(p-h)} }{h R_b} \right) - \sqrt{p(p-R_b)^2(p-h)} \bigg).
$$
We finally obtain
\Beq
\label{Chapriskmeasures_Aire_Inter}
A_{int}(h)=
\left\{
\begin{array}{l l}
2 \bigg(R_b^2 \arcsin \left( \frac{ \sqrt{4R_b^2-h^2}}{ 2 R_b}  \right) - \frac{h}{4} \sqrt{4 R_b^2-h^2} \bigg) & \mbox{ if } \  h \leq 2 R_b,  \\
0 & \mbox{ if } \  h > 2 R_b.
\end{array}
\right.
\Eeq
The combination of \eqref{Bivariate_Distri_Function_Tube} and \eqref{Chapriskmeasures_Aire_Inter} yields the result.
\end{proof}

\subsection{For Corollary \ref{Maxstable_Tubes}}

\begin{proof}
By definition of the extremal coefficient function, for all $\Mb{x}_1, \Mb{x}_2 \in \mathds{R}^2$, $\Theta(\Mb{x}_1, \Mb{x}_2) =- \log (\mathds{P}(Z(\Mb{x}_1) \leq u, Z(\Mb{x}_2) \leq u) \ u.$
Hence, using \eqref{Bivariate_Distri_Function_Tube}, we obtain
\Beq
\label{Chapriskmeasures_Extremal_Coefficient_Tube_Firstexpression}
\Theta(\Mb{x}_1, \Mb{x}_2) = u \left( \frac{h_b}{u} \left[ \pi R_b^2-A_{int}(h) \right] + \frac{1}{u}  \right)=2-h_b \  A_{int}(h).
\Eeq
The combination of \eqref{Chapriskmeasures_Aire_Inter} and \eqref{Chapriskmeasures_Extremal_Coefficient_Tube_Firstexpression} gives the result.
\end{proof}

\subsection{For Proposition \ref{Proposition_Variance_Tubes}}

\begin{proof}
We consider the case of $A$ being a disk; the proof is exactly the same in the case of a square.
Using \eqref{Risk_Measure_Disk}, we have
\begin{align*}
\mathcal{R}_2(\lambda A) 
&= - \exp \left( -\frac{2}{u} \right)
+ \int_{0}^{2 R_b} f_d(h,R) \ \exp \left( -\frac{\Theta(\lambda h)}{u} \right) \  dh 
\\& \quad \ + \int_{2 R_b}^{2 R} f_d(h,R) \ \exp \left( -\frac{\Theta(\lambda h)}{u} \right) \  dh.
\end{align*}
Moreover, for $\lambda \geq 1$ and $h \geq 2 R_b$, we have $\Theta(\lambda h)=\Theta(h)=2$. Hence, for all $\lambda \geq 1$,
\begin{align*}
& \quad \ \mathcal{R}_2 (\lambda A) 
\\& = \int_{0}^{2 R_b} f_d(h,R) \ \exp \left( -\frac{\Theta(\lambda h)}{u} \right) \  dh - \exp \left( -\frac{2}{u} \right) + \int_{0}^{2 R} f_d(h,R) \ \exp \left( -\frac{\Theta(h)}{u} \right) \  dh 
\\& \ \ \ - \int_{0}^{2 R_b} f_d(h,R) \ \exp \left( -\frac{\Theta(h)}{u} \right) \  dh
\\& = \mathcal{R}_2(A) + \int_{0}^{2 R_b} f_d(h,R) \left[  \exp \left( -\frac{\Theta(\lambda h)}{u} \right)-\exp \left( -\frac{\Theta(h)}{u} \right) \right] dh.
\end{align*}
When $R_b$ tends to 0, the second term vanishes, giving $\lim_{R_b \to 0} \mathcal{R}_2(\lambda A)=\mathcal{R}_2(A)$. Moreover, if $R_b=0$, for all $h \geq 0$, for all $\lambda \geq 0, \Theta(\lambda h)=2$. Thus, \eqref{Risk_Measure_Disk} gives that for all $\lambda \geq 0, \mathcal{R}_2(\lambda A)=0$.
\end{proof}

\subsection{For Theorem \ref{Variance_Properties}}

\begin{proof} 
1. If $Z$ is stationary, it is also true for the process $C_P=\mathds{I}_{\{ Z>u  \}}$. All processes mentioned are stationary. See e.g. \cite{schlather2002models}, Theorem 2, for the Schlather process and \cite{kabluchko2009stationary}, Theorem 2, for the Brown-Resnick process; the tube process is obviously stationary as an instance of M2 process.

\medskip

\noindent 2. Let us consider two regions $A_1$ and $A_2$ being both a disk or a square and such that $A_1 \subset A_2$. Due to spatial invariance under translation, the region $A_2$ can be translated to region $A_2'$, where $A_2'$ corresponds to the region obtained by an homothety of $A_1$, with center the center of $A_1$ and ratio denoted $\lambda \geq 1$. Thus, $\mathcal{R}_2(A_2)=\mathcal{R}_2(A_2')=\mathcal{R}_2(\lambda A_1)$. Moreover, for all $h>0$, the function $\lambda \mapsto \theta(\lambda h)$ is non-decreasing. Thus, \eqref{Risk_Measure_Disk} and \eqref{Risk_Measure_Square} give that $\mathcal{R}_2(\lambda A)$ is a non-increasing function of $\lambda$. Hence, we have $\mathcal{R}_2(\lambda A_1) \leq \mathcal{R}_2(A_1)$, giving $\mathcal{R}_2(A_2) \leq \mathcal{R}_2(A_1)$. Hence, there is spatial anti-monotonicity, directly implying spatial sub-additivity. All processes mentioned are stationary and have an isotropic and non-decreasing extremal coefficient function. For the Brown-Resnick model, we have, in the general case, $\Theta(\Mb{x}, \Mb{y})=2 \Phi \left( \sqrt{\frac{\gamma(\Mb{x}-\Mb{y})}{2}} \right)$, $\Mb{x}, \Mb{y} \in \mathds{R}^2$ \citep[see e.g.][]{davison2012statistical}.

\medskip

\noindent 3.
First, we show that
\Beq
\label{Eq_Equivalence_Conditions}
0 < \int_{\mathds{R}^2} (2-\Theta(\Mb{x})) \  d\Mb{x} < \infty   \Leftrightarrow 0 < \sigma^2 < \infty.
\Eeq
Since $u>0$ and, for all $\Mb{x} \in \mathds{R}^2$, $\Theta(\Mb{x}) \in [1,2]$, we have, using the mean value theorem, that
$$
0 \leq \exp \left( -\frac{\Theta(\Mb{x})}{u} \right)-\exp \left( -\frac{2}{u} \right) \leq \frac{2-\Theta(\Mb{x})}{u}.
$$
Hence, we immediately obtain that
\Beq
\label{Eq_Implication_1}
\int_{\mathds{R}^2} (2-\Theta(\Mb{x})) \  d\Mb{x} < \infty \Rightarrow \sigma^2 < \infty.
\Eeq
Moreover, we have that
\begin{align*}
\exp \left( -\frac{\Theta(\Mb{x})}{u} \right) - \exp \left( -\frac{2}{u} \right)& = \exp \left( -\frac{2}{u} \right) \left( \exp \left( \frac{2-\Theta(\Mb{x})}{u} \right) -1 \right) 
\\& \geq  \exp \left( -\frac{2}{u} \right) \frac{2-\Theta(\Mb{x})}{u},
\end{align*}
which directly yields that
\Beq
\label{Eq_Implication_2}
\sigma^2 < \infty \Rightarrow \int_{\mathds{R}^2} (2-\Theta(\Mb{x})) \  d\Mb{x} < \infty.
\Eeq
Moreover, since for all $\Mb{x} \in \mathds{R}^2$, $2-\Theta(\Mb{x}) \geq 0$,  $\displaystyle \int_{\mathds{R}^2} (2-\Theta(\Mb{x})) \  d\Mb{x} >0$ iff $2-\Theta(\Mb{x})>0$ for all $\Mb{x} \in \mathcal{D}$, where $\mathcal{D} \subset \mathds{R}^2$ and satisfies $| \mathcal{D} |>0$. This directly gives
\Beq
\label{Eq_Equivalence_2}
\int_{\mathds{R}^2} (2-\Theta(\Mb{x})) \  d\Mb{x} >0 \Leftrightarrow \sigma^2>0.
\Eeq
The combination of \eqref{Eq_Implication_1}, \eqref{Eq_Implication_2} and \eqref{Eq_Equivalence_2} gives \eqref{Eq_Equivalence_Conditions}. Thus, \eqref{Eq_Condition_Thm_3} and \eqref{Eq_Equivalence_Conditions} yield that $0 < \sigma^2 < \infty$. Note that, although only the ``$\Rightarrow$'' implication of  \eqref{Eq_Equivalence_Conditions} is directly useful in this proof, we have also proven the ``$\Leftarrow$'' implication in order to show that \eqref{Eq_Condition_Thm_3} is not more restrictive than the condition $0 < \sigma^2 < \infty$. 

\medskip

The remainder of the proof is partly inspired by the proof of Proposition 7.3 in \cite{dedecker2009weak}.
Let $A$ be a convex in $\mathcal{A}$. For $\lambda>0$, we denote $A_{\lambda}=\lambda A$. For a set $F \subset \mathds{R}^2$, let $\partial F$ and $P_F$ be the boundary and the perimeter of $F$, respectively. 
Let 
\Beq
\label{Eq_Equation_vlambda}
v_\lambda=\lambda^{\alpha},
\Eeq 
for some $\alpha \in (0,1)$. Let also
$\partial_{v_{\lambda}} A_\lambda= \{ \Mb{x} \in A_\lambda: \mbox{dist} (\Mb{x}, \partial A_\lambda) \leq v_\lambda \}$ and $A_\lambda^{(v_\lambda)} = A_\lambda \backslash \partial_{v_{\lambda}} A_\lambda$.
Since $A_{\lambda}$ is compact, convex and has a positive Lebesgue measure, we easily obtain, using Theorem 2 (note that this result is related to Steiner's formula) in \cite{swanson2011band}, that 
$| \partial_{v_{\lambda}} A_\lambda | = v_{\lambda} P_{A_\lambda} - \pi v_{\lambda}^2=\lambda v_{\lambda} P_A -\pi v_{\lambda}^2 \leq \lambda v_{\lambda} P_A.$
Hence, 
$$\frac{| \partial_{v_{\lambda}} A_\lambda |}{\lambda^2 |A|} \leq \frac{v_{\lambda} P_A}{\lambda |A|},$$
which gives, using \eqref{Eq_Equation_vlambda}, that
\Beq
\label{Eq_Lim_Partial}
\lim_{\lambda \to \infty} \frac{| \partial_{v_{\lambda}} A_\lambda |}{\lambda^2 |A|}=0.
\Eeq
By definition of $A_\lambda^{(v_\lambda)}$, we obtain that $|A_\lambda^{(v_\lambda)}|=\lambda^2 |A|- | \partial_{v_{\lambda}} A_\lambda |$, which gives, using \eqref{Eq_Lim_Partial}, that
\Beq
\label{Eq_Lim_Comp}
\lim_{\lambda \to \infty} \frac{|A_\lambda^{(v_\lambda)}|}{\lambda^2 |A|}=1.
\Eeq

Using \eqref{Chapriskmeasures_Eq_Stationary_Risk_Measure_Generalcase} and
the fact that, for all $\Mb{x}, \Mb{y} \in \mathds{R}^2$, $\Theta(\Mb{x}, \Mb{y})$ can be written $\Theta(\Mb{x}-\Mb{y})$, we have
\Beq
\label{Eq_R2_Theta_Stationary}
\mathcal{R}_2(\lambda A)= \frac{1}{\lambda^4 |A|^2} \int_{\lambda A}  \int_{\lambda A} \left[ \exp\left( -\frac{\Theta(\Mb{x}-\Mb{y})}{u} \right) - \exp \left( -\frac{2}{u} \right) \right] \  d\Mb{x}\  d\Mb{y}.
\Eeq
We introduce
$$ k(\Mb{x})= \exp \left( -\frac{\Theta(\Mb{x})}{u} \right)  -  \exp \left( -\frac{2}{u} \right), \quad \Mb{x} \in \mathds{R}^2$$
and $$T_\lambda = \frac{1}{\lambda^2 |A|} \int_{A_\lambda} \int_{A_\lambda} k(\Mb{x}-\Mb{y})\  d\Mb{x} \ d\Mb{y}, \quad \lambda >0.$$
Using \eqref{Eq_R2_Theta_Stationary} and \eqref{Eq_Def_Sigma}, these notations lead to
\Beq
\label{Eq_Notations}
\mathcal{R}_2(\lambda A)= \frac{1}{\lambda^2 |A|} T_{\lambda} \quad \mbox{and} \quad \sigma^2= \int_{\mathds{R}^2} k(\Mb{x}) \ d\Mb{x}.
\Eeq
Note that, for all $\Mb{x} \in \mathds{R}^2$, $\Theta( \Mb{x} ) \leq 2$, meaning that $k$ is non-negative. We first show that $\lim_{\lambda \to \infty} T_\lambda=\sigma^2$. We have, for all $\lambda>0$, that
\Beq
\label{Eq_T_Eq1}
T_\lambda = T_{1,\lambda} + T_{2,\lambda} + T_{3,\lambda}, 
\Eeq
where
\begin{align*}
& T_{1, \lambda} = \frac{1}{\lambda^2 |A|} \int_{A_\lambda^{(v_\lambda)}} \int_{\Mb{y} \in A_\lambda: \| \Mb{x} - \Mb{y} \| \geq v_\lambda } k(\Mb{x}-\Mb{y}) \  d\Mb{x} \  d\Mb{y},
\\& T_{2, \lambda} = \frac{1}{\lambda^2 |A|} \int_{A_\lambda^{(v_\lambda)}} \int_{\Mb{y} \in A_\lambda: \| \Mb{x} - \Mb{y} \| < v_\lambda } k(\Mb{x}-\Mb{y}) \  d\Mb{x} \  d\Mb{y},
\\& T_{3, \lambda} = \frac{1}{\lambda^2 |A|} \int_{ \partial_{v_{\lambda}} A_{\lambda}} \int_{A_\lambda} k(\Mb{x}-\Mb{y}) \  d\Mb{x} \  d\Mb{y}.
\end{align*}

\noindent \textbf{Study of $T_{1,\lambda}$}

Using the facts that $k$ is non-negative and that $|A_\lambda^{(v_\lambda)}| \leq \lambda^2 |A|$, we obtain, for all $\lambda > 0$, that
\begin{align*}
| T_{1,\lambda} | & \leq \frac{1}{\lambda^2 |A|} \int_{A_\lambda^{(v_\lambda)}} \int_{\Mb{y} \in \mathds{R}^2: \| \Mb{x} - \Mb{y} \| \geq v_\lambda } k(\Mb{x}-\Mb{y}) \  d\Mb{x} \  d\Mb{y} 
\\& = \frac{1}{\lambda^2 |A|} \int_{A_\lambda^{(v_\lambda)}} \left( \int_{\Mb{z} \in \mathds{R}^2: \| \Mb{z} \| \geq v_\lambda } k(\Mb{z}) \  d\Mb{z} \right) \ d\Mb{x}
\\& \leq \int_{\Mb{z} \in \mathds{R}^2: \| \Mb{z} \| \geq v_\lambda } k(\Mb{z}) \  d\Mb{z}.
\end{align*}
Thus, combining the fact that $\lim_{\lambda \to \infty} v_{\lambda}= \infty$, the second part of \eqref{Eq_Notations} and the fact that $\sigma^2$ is finite yields that
\Beq
\label{Eq_T1_Eq3}
\lim_{\lambda \to \infty} T_{1,\lambda}=0.
\Eeq

\noindent \textbf{Study of $T_{2,\lambda}$}

If $\Mb{x} \in A_\lambda^{(v_\lambda)}$ and $\Mb{y} \notin A_\lambda$, we necessarily have that $\| \Mb{x} - \Mb{y} \| \geq v_\lambda$.
Thus, we have, for all $\lambda>0$, that
\begin{align*}
T_{2,\lambda} &= \frac{1}{\lambda^2 |A|} \int_{A_\lambda^{(v_\lambda)}} \int_{\Mb{y} \in \mathds{R}^2: \| \Mb{x} - \Mb{y} \| < v_\lambda } k(\Mb{x}-\Mb{y}) \  d\Mb{x} \  d\Mb{y}
\\&= \frac{1}{\lambda^2 |A|} \int_{A_\lambda^{(v_\lambda)}} \left( \int_{\Mb{z} \in \mathds{R}^2: \| \Mb{z} \| < v_\lambda } k(\Mb{z}) \ d\Mb{z} \right) \ d\Mb{x}
\\& = \frac{|A_\lambda^{(v_\lambda)}|}{\lambda^2 |A|} \int_{\Mb{z} \in \mathds{R}^2: \| \Mb{z} \| < v_\lambda } k(\Mb{z}) \ d\Mb{z}.
\end{align*}
Hence, we obtain, using \eqref{Eq_Lim_Comp}, the fact that $\lim_{\lambda \to \infty} v_{\lambda}=\infty$ and the second part of \eqref{Eq_Notations}, that 
\Beq
\label{Eq_T2}
\lim_{\lambda \to \infty} T_{2,\lambda} = \sigma^2.
\Eeq

\noindent \textbf{Study of $T_{3,\lambda}$}

We have, using the fact that $k$ is non-negative and the second part of \eqref{Eq_Notations}, that, for all $\lambda>0$,
\begin{align*}
| T_{3,\lambda} | \leq \frac{1}{\lambda^2 |A|} \int_{ \partial_{v_{\lambda}} A_{\lambda}} \int_{\mathds{R}^2} k(\Mb{x}-\Mb{y}) \  d\Mb{x} \  d\Mb{y} = \frac{1}{\lambda^2 |A|} \int_{ \partial_{v_{\lambda}} A_{\lambda}}  \left( \int_{\mathds{R}^2} k(\Mb{z}) \  d\Mb{z} \right) \  d\Mb{x} = \frac{| \partial_{v_{\lambda}} A_{\lambda} |}{\lambda^2 |A|} \sigma^2.
\end{align*}
Therefore, using \eqref{Eq_Lim_Partial}, we obtain
\Beq
\label{Eq_T3}
\lim_{\lambda \to \infty} T_{3,\lambda} = 0.
\Eeq

\noindent Finally, the combination of \eqref{Eq_T_Eq1}, \eqref{Eq_T1_Eq3}, \eqref{Eq_T2} and \eqref{Eq_T3} gives that
$$\lim_{\lambda \to \infty} T_\lambda = \sigma^2.$$

Thus, using the first part of \eqref{Eq_Notations}, we obtain
$$ \lambda^2 |A| \mathcal{R}_2 (\lambda A) \underset{\lambda \to \infty}{=} \sigma^2 + o(1)
\quad \mbox{ i.e. } \quad \mathcal{R}_2(\lambda A) \underset{\lambda \to \infty}{=} \dfrac{\sigma^2}{\lambda^2 |A|}+o \left( \dfrac{1}{\lambda^2} \right).$$ Consequently, setting $K_1=0$ and $K_2=\dfrac{\sigma^2}{|A|}$, the axiom of asymptotic spatial homogeneity of order $-2$ is satisfied since $\sigma^2 >0$.

\medskip

Let us now prove that the Smith process, the Brown-Resnick process with a semivariogram satisfying $\gamma( \Mb{h} ) = \eta \| \Mb{h} \|^a$ where $\eta>0$ and $a \in (0, 2]$, and the tube process satisfy the required conditions. First, they are all simple and stationary. Thus, for all $\Mb{x}, \Mb{y} \in \mathds{R}^2$, $\Theta(\Mb{x}, \Mb{y})$ only depends on $\Mb{x}-\Mb{y}$. Moreover, for all $\Mb{x} \in \mathds{R}^2$, $2-\Theta(\Mb{x})>0$ on a set with positive Lebesgue measure. Hence, $\displaystyle \int_{\mathds{R}^2} (2-\Theta(\Mb{x})) \  d\Mb{x} >0$. \\
\noindent We now show that $\displaystyle \int_{\mathds{R}^2} (2-\Theta(\Mb{x})) \  d\Mb{x} < \infty$. In the case of the Smith process, recall that 
\Beq
\label{Eq_Smith_Coeff_Extremal}
\Theta(\Mb{x})=2 \Phi \left( \frac{\| \Mb{x} \|_{\Sigma} }{2} \right), \quad \Mb{x} \in \mathds{R}^2,
\Eeq
where $\| \Mb{x} \|_{\Sigma}=\sqrt{\Mb{x}^{'} \Sigma^{-1} \Mb{x}}$. Note that $\|.\|_{\Sigma}$ is the norm associated with the inner product induced by the matrix $\Sigma^{-1}$.
Mill's ratio gives us that the survival distribution function $\bar{\Phi}$ of a standard Gaussian random variable satisfies
$\bar{\Phi}(h) \underset{h \to \infty}{\sim} \dfrac{\exp \left( -\frac{h^2}{2} \right)}{\sqrt{2 \pi} h}$.
Hence, using \eqref{Eq_Smith_Coeff_Extremal}, we obtain
\Beq
\label{Eq_Equivalence_f}
2-\Theta(\Mb{x}) = 2 \left[ 1-\Phi \left( \frac{\| \Mb{x} \|_{\Sigma}}{2} \right) \right] \underset{\| \Mb{x} \| \to \infty}{\sim} \frac{ 4 \exp \left( -\frac{\| \Mb{x} \|_{\Sigma}^2}{8} \right) }{ \sqrt{2 \pi} \| \Mb{x} \|_{\Sigma} }.
\Eeq
Since all norms in $\mathds{R}^2$ are equivalent, there exists $B_1>0$ such that, for all $\Mb{x} \in \mathds{R}^2$, $\| \Mb{x} \|_{\Sigma} \geq B_1 \| \Mb{x} \|$. Thus, for all $\Mb{x} \in \mathds{R}^2$, 
\Beq
\label{Eq_Integrability_Intermediate_Term}
0 \leq \frac{ 4 \exp \left( -\frac{\| \Mb{x} \|_{\Sigma}^2}{8} \right) }{ \sqrt{2 \pi} \| \Mb{x} \|_{\Sigma}} \leq \frac{ 4 \exp \left( -\frac{B_1^2 \| \Mb{x} \|^2}{8} \right) }{ \sqrt{2 \pi} B_1 \| \Mb{x} \|}.
\Eeq
Using the fact that the right-hand term of \eqref{Eq_Integrability_Intermediate_Term} is integrable, the combination of \eqref{Eq_Equivalence_f} and \eqref{Eq_Integrability_Intermediate_Term} gives the integrability of $2-\Theta$. \\
\noindent In the case of the Brown-Resnick process with a semivariogram satisfying $\gamma( \Mb{h} ) = \eta \| \Mb{h} \|^a$ where $\eta>0$ and $a \in (0, 2]$, the integrability is obtained in exactly the same way since $\Theta(\Mb{x})=2 \Phi \left( \sqrt{\frac{\gamma(\Mb{x})}{2}} \right)$, $\Mb{x} \in \mathds{R}^2$. \\
\noindent In the case of the tube model, we have $\Theta( \Mb{x} )=2$ for $\| \Mb{x} \| \geq 2R_b$. Therefore, the function $2-\Theta$ has a compact support and is integrable. \\
\noindent Finally, in all cases mentioned, we have $\displaystyle \int_{\mathds{R}^2} (2-\Theta(\Mb{x})) \  d\Mb{x} < \infty$. 
\end{proof}

\subsection{For Theorem \ref{TCL_Temperature}} 

\begin{proof}
A random field $\{X(\Mb{x})\}_{\Mb{x} \in \mathds{R}^2}$ is called associated  if 
$\mbox{Cov} ( f(X_I), g(X_I) ) \geq 0$
for any finite subset $I \subset \mathds{R}^2$ and for any bounded coordinatewise non-decreasing  functions $f:\mathds{R}^{\mbox{card}(I)} \mapsto \mathds{R}$ and $g:\mathds{R}^{\mbox{card}(I)} \mapsto \mathds{R}$, where $X_I=\{X(\Mb{x}): \Mb{x} \in I\}$ and card stands for cardinality. A random field $\{X(\Mb{x})\}_{\Mb{x} \in \mathds{R}^2}$ is called positively associated if $\mbox{Cov}( f(X_I), g(X_J) ) \geq 0$
for all disjoint finite subsets $I,J \subset \mathds{R}^2$ and for any bounded coordinatewise non-decreasing  functions $f:\mathds{R}^{\mbox{card}(I)} \mapsto \mathds{R}$ and $g:\mathds{R}^{\mbox{card}(J)} \mapsto \mathds{R}$. If $X$ is associated, it is clear that $X$ is positively associated. Proposition 5.29 in \cite{resnickextreme} implies that max-stable processes are associated and therefore positively associated. \\
\noindent Moreover, $Z$ is stationary and measurable and thus continuous in probability. Since $Z$ has standard Fr\'echet margins, we obviously have that $\mathds{P}(Z(\Mb{0})=u)=0$ for all $u>0$. \\
\noindent For $r>0$, we denote $F^{+r}=\{ \Mb{x} \in \mathds{R}^2: \mbox{dist}(\Mb{x}, F) \leq r \}$,
where $\mbox{dist}$ is the Euclidean distance. A Van Hove sequence in $\mathds{R}^2$ is a sequence $( F_n )_{n \in \mathds{N}}$ of measurable subsets of $\mathds{R}^2$ satisfying 
$F_n \uparrow \mathds{R}^2$, $\lim_{n \to \infty} |F_n|=\infty$, and $\lim_{n \to \infty} \dfrac{| (\partial F_n)^{+r} |}{|F_n|} =0, \mbox{ for all } r>0.$
For any convex $A \in \mathcal{A}$ and any positive non-decreasing sequence $(\lambda_n)_{n \in \mathds{N}}$ such that $\lim_{n \to \infty} \lambda_n=\infty$, the sequence $(\lambda_n A)_{n \in \mathds{N}}$ is a Van Hove sequence. Indeed, for all $n \in \mathds{N}$, $\lambda_n A \subset \lambda_{n+1}A$ and $\bigcup_{n=1}^{\infty} \lambda_n A=\mathds{R}^2$, showing that the first condition is satisfied. The second one is trivial and we now show the last one. For all $\lambda>0$, since $A_{\lambda}$ is compact, convex and has a positive Lebesgue measure, we easily obtain, using Theorem 2 in \cite{swanson2011band}, that $| (\partial A_{\lambda})^{+r} | = 2r P_{A_{\lambda}} =2 r \lambda P_A$.
This directly gives
$$
\lim_{\lambda \to \infty} \frac{|(\partial A_\lambda)^{+r}|}{\lambda^2 |A|}=0,
$$
which immediately implies the last condition in the definition of a Van Hove sequence.\\
\noindent Furthermore, it turns out that
\begin{align*}
& \quad \ \int_{\mathds{R}^2} \mbox{Cov} \left( \mathds{I}_{\{Z(\Mb{0})>u\}}, \mathds{I}_{\{Z(\Mb{x})>u\}} \right) d\Mb{x}
\\&= \int_{\mathds{R}^2} \left( \mathds{P}(Z(\Mb{0})>u, Z(\Mb{x})>u) - \left[ 1-\exp \left( -\frac{1}{u} \right) \right]^2 \right) d\Mb{x}
\\& = \int_{\mathds{R}^2} \left[ \exp \left( -\frac{\Theta(\Mb{x})}{u} \right) - \exp \left( -\frac{2}{u} \right) \right] d\Mb{x}=\sigma^2. 
\end{align*} \\
\noindent The combination of \eqref{Eq_Condition_Thm_4} and \eqref{Eq_Implication_1} implies that $\sigma^2 < \infty$. Moreover, $\sigma^2 >0$ since $Z$ is measurable. Indeed, the case $\sigma^2=0$ corresponds to the situation where the random variables $Z(\Mb{x})$ are independent for almost all $\Mb{x} \in \mathds{R}^2$, which corresponds to processes which are not measurable. Therefore, we have that $\sigma^2 \in (0, \infty)$ and, by applying Theorem 7 in \cite{spodarev2014limit}, we obtain, for all positive non-decreasing sequence $(\lambda_n)_{n \in \mathds{N}}$ such that $\lim_{n \to \infty} \lambda_n=\infty$, that, for all convex $A \in \mathcal{A}$,
\Beq
\label{Eq_TCL_Lambda_n}
\lambda_n \sqrt{|A|} \left( L_N(\lambda_n A) - m \right) \overset{d}{\to} \mathcal{N}(0, \sigma^2), \quad \mbox{i.e.} \quad \lambda_n \left( L_N(\lambda_n A) - m \right) \overset{d}{\to} \mathcal{N} \left( 0, \frac{\sigma^2}{|A|} \right),
\Eeq
for $n \to \infty$.
Now, let $l: \mathds{R} \to \mathds{R}$ and $B_2 \in \mathds{R}$. If, for any positive non-decreasing sequence $(\lambda_n)_{n \in \mathds{N}}$ such that $\lim_{n \to \infty} \lambda_n=\infty$, we have that $\lim_{n \to \infty} l(\lambda_n)=B_2$, then $\lim_{\lambda \to \infty} l(\lambda) = B_2$. To prove this, let us assume that there exists a sequence $(\eta_n)_{n \in \mathds{N}}$ such that $\lim_{n \to \infty} \eta_n=\infty$ but $(l(\eta_n))_{n \in \mathds{N}}$ does not converge to $B_2$. This means that there exists $\epsilon>0$ such that, for all $N \in \mathds{N}$, there exists $n \geq N$ such that $|l(\eta_n)-B_2|>\epsilon$. We introduce the sequence $(N_n)_{n \in \mathds{N}}$ defined by $N_1=\inf \{ n \in \mathds{N} : \eta_n >0 \mbox{ and } |l(\eta_n)-B_2|> \epsilon \}$ and, for $n \in \mathds{N}$, $N_{n+1} = \inf \{ n \in \mathds{N} \backslash \{ 1, \dots, N_n \}: |l(\eta_n)-B_2|> \epsilon \mbox{ and } \eta_n \geq \eta_{N_n} \}$. It is clear that $N_1$ and $N_{n+1}$ exist since $\lim_{n \to \infty} \eta_n=\infty$. By definition, the sequence $(\eta_{N_n})_{n \in \mathds{N}}$ is positive, non-decreasing, satisfies $\lim_{n \to \infty} \eta_{N_n}=\infty$ (as a subsequence of $(\eta_n)_{n \in \mathds{N}}$) but does not converge to $B_2$. This leads to a contradiction and shows the result. 
Hence, finally, we obtain \eqref{Eq_Convergence_TCL} using \eqref{Eq_TCL_Lambda_n}.

\medskip

Finally, the Smith process, the Brown-Resnick process with a semivariogram satisfying $\gamma( \Mb{h} ) = \eta \| \Mb{h} \|^a$ where $\eta>0$ and $a \in (0, 2]$, and the tube process are all stationary and measurable. Moreover, in all these cases, as shown in the proof of Theorem \ref{Variance_Properties}, \eqref{Eq_Condition_Thm_4} is satisfied. Thus, all the required conditions are satisfied.
\end{proof}

\subsection{For Theorem \ref{Th_Asymptotic_Homogeneity_VaR}}

\begin{proof}
1. The argument is the same as in the proof of Theorem \ref{Variance_Properties}, Bullet 1.

\medskip

\noindent 2.
By Theorem \ref{TCL_Temperature} and Proposition 0.1 in \cite{resnickextreme}, we have that
$$ \mbox{VaR}_{\alpha} ( \lambda (L_N(\lambda A))-m) ) = \lambda ( \mbox{VaR}_{\alpha}(L_N(\lambda A)) - m ) \to \frac{\sigma \ q_{\alpha}}{\sqrt{|A|}}, \mbox{ for } \lambda \to \infty.$$
Thus,
$$ \lambda ( \mbox{VaR}_{\alpha}(L_N(\lambda A)) - m ) \underset{\lambda \to \infty}{=} \frac{\sigma \ q_{\alpha}}{\sqrt{|A|}} + o(1) \quad \mbox{ i.e. } \quad
\mathcal{R}_{3, \alpha}(\lambda A) \underset{\lambda \to \infty}{=} m+ \frac{\sigma \  q_{\alpha}}{\lambda \sqrt{|A|}} +o \left( \frac{1}{\lambda} \right).
$$
Hence, setting $K_1=m$ and $K_2=\dfrac{\sigma \ q_{\alpha}}{\sqrt{|A|}}$, the axiom of asymptotic spatial homogeneity of order $-1$ is satisfied since $\sigma>0$ and $q_{\alpha} \neq 0$ (because $\alpha \neq 0.5$). 

\medskip

Finally, regarding the examples of processes proposed, we know from Theorem \ref{TCL_Temperature} that the required conditions are satisfied.
\end{proof}

\newpage
\bibliographystyle{apalike}
\bibliography{Biblio}

\begin{thebibliography}{}

\bibitem[Abrahamsen, 1997]{abrahamsen1997review}
Abrahamsen, P. (1997).
\newblock A review of {G}aussian random fields and correlation functions.
\newblock Technical Report 917, Norwegian Computing Center.

\bibitem[Acciaio and Penner, 2011]{acciaio2011dynamic}
Acciaio, B. and Penner, I. (2011).
\newblock Dynamic risk measures.
\newblock In {\em Advanced Mathematical Methods for Finance}, pages 1--34.
  Springer.

\bibitem[Adler et~al., 2013]{adler2013high}
Adler, R.~J., Samorodnitsky, G., and Taylor, J.~E. (2013).
\newblock High level excursion set geometry for non-{G}aussian infinitely
  divisible random fields.
\newblock {\em The Annals of Probability}, 41(1):134--169.

\bibitem[Akaike, 1974]{akaike1974new}
Akaike, H. (1974).
\newblock A new look at the statistical model identification.
\newblock {\em IEEE Transactions on Automatic Control}, 19(6):716--723.

\bibitem[Artzner et~al., 1999]{artzner1999coherent}
Artzner, P., Delbaen, F., Eber, J.-M., and Heath, D. (1999).
\newblock Coherent measures of risk.
\newblock {\em Mathematical Finance}, 9(3):203--228.

\bibitem[Bevere and Mueller, 2014]{SwissRe}
Bevere, L. and Mueller, L. (2014).
\newblock Natural catastrophes and man-made disasters in 2013: large losses
  from floods and hail; {H}aiyan hits the {P}hilippines.
\newblock {\em Sigma Swiss Re}, 2014(1).

\bibitem[Brown and Resnick, 1977]{brown1977extreme}
Brown, B.~M. and Resnick, S.~I. (1977).
\newblock Extreme values of independent stochastic processes.
\newblock {\em Journal of Applied Probability}, 14(4):732--739.

\bibitem[Cousin and Di~Bernardino, 2013]{cousin2013multivariate}
Cousin, A. and Di~Bernardino, E. (2013).
\newblock On multivariate extensions of {V}alue-at-{R}isk.
\newblock {\em Journal of Multivariate Analysis}, 119:32--46.

\bibitem[Davison, 2003]{davison2003statistical}
Davison, A.~C. (2003).
\newblock {\em Statistical {M}odels}.
\newblock Cambridge University Press.

\bibitem[Davison et~al., 2012]{davison2012statistical}
Davison, A.~C., Padoan, S.~A., and Ribatet, M. (2012).
\newblock Statistical modeling of spatial extremes.
\newblock {\em Statistical Science}, 27(2):161--186.

\bibitem[de~Haan, 1984]{haan1984spectral}
de~Haan, L. (1984).
\newblock A spectral representation for max-stable processes.
\newblock {\em The Annals of Probability}, 12(4):1194--1204.

\bibitem[de~Haan and Ferreira, 2006]{de2007extreme}
de~Haan, L. and Ferreira, A. (2006).
\newblock {\em Extreme {V}alue {T}heory: {A}n {I}ntroduction}.
\newblock Springer.

\bibitem[de~Haan and Pickands, 1986]{de1986stationary}
de~Haan, L. and Pickands, J. (1986).
\newblock Stationary min-stable stochastic processes.
\newblock {\em Probability Theory and Related Fields}, 72(4):477--492.

\bibitem[Dedecker et~al., 2007]{dedecker2009weak}
Dedecker, J., Doukhan, P., Lang, G., Le{\'o}n~R, J.~R., Louhichi, S., and
  Prieur, C. (2007).
\newblock {\em Weak Dependence: With Examples and Applications}.
\newblock Springer.

\bibitem[Dombry, 2012]{DombryHDR2012}
Dombry, C. (2012).
\newblock {\em Th\'eorie spatiale des extr\^emes et propri\'et\'es des
  processus max-stables}.
\newblock Habilitation \`a diriger des recherches ({HDR}), Universit\'e de
  Poitiers.

\bibitem[Embrechts et~al., 2005]{embrechts2005strategic}
Embrechts, P., Kaufmann, R., and Patie, P. (2005).
\newblock Strategic long-term financial risks: Single risk factors.
\newblock {\em Computational Optimization and Applications}, 32(1):61--90.

\bibitem[Embrechts and Puccetti, 2006]{embrechts2006bounds}
Embrechts, P. and Puccetti, G. (2006).
\newblock Bounds for functions of multivariate risks.
\newblock {\em Journal of Multivariate Analysis}, 97(2):526--547.

\bibitem[F{\"o}llmer, 2014]{follmer2014spatial}
F{\"o}llmer, H. (2014).
\newblock Spatial risk measures and their local specification: the locally
  law-invariant case.
\newblock {\em Statistics \& Risk Modeling}, 31(1):79--101.

\bibitem[F{\"o}llmer and Kl{\"u}ppelberg, 2014]{follmerspatial}
F{\"o}llmer, H. and Kl{\"u}ppelberg, C. (2014).
\newblock Spatial risk measures: Local specification and boundary risk.
\newblock In {\em Stochastic Analysis and Applications 2014}, pages 307--326.
  Springer.

\bibitem[F{\"o}llmer and Schied, 2002]{follmer2002convex}
F{\"o}llmer, H. and Schied, A. (2002).
\newblock Convex measures of risk and trading constraints.
\newblock {\em Finance and {S}tochastics}, 6(4):429--447.

\bibitem[F\"ollmer and Schied, 2004]{Follmer2004}
F\"ollmer, H. and Schied, A. (2004).
\newblock {\em Stochastic Finance: An Introduction in Discrete Time}.
\newblock de Gruyter.

\bibitem[Frittelli and Rosazza~Gianin, 2002]{frittelli2002putting}
Frittelli, M. and Rosazza~Gianin, E. (2002).
\newblock Putting order in risk measures.
\newblock {\em Journal of Banking \& Finance}, 26(7):1473--1486.

\bibitem[Gorodetskii, 1987]{gorodetskii1987moment}
Gorodetskii, V.~V. (1987).
\newblock Moment inequalities and the central limit theorem for integrals of
  random fields with mixing.
\newblock {\em Journal of Soviet Mathematics}, 36(4):461--467.

\bibitem[Guidolin and Timmermann, 2006]{guidolin2006term}
Guidolin, M. and Timmermann, A. (2006).
\newblock Term structure of risk under alternative econometric specifications.
\newblock {\em Journal of Econometrics}, 131(1):285--308.

\bibitem[Kabluchko, 2009]{kabluchko2009spectral}
Kabluchko, Z. (2009).
\newblock Spectral representations of sum- and max-stable processes.
\newblock {\em Extremes}, 12(4):401--424.

\bibitem[Kabluchko and Schlather, 2010]{kabluchko2010ergodic}
Kabluchko, Z. and Schlather, M. (2010).
\newblock Ergodic properties of max-infinitely divisible processes.
\newblock {\em Stochastic Processes and their Applications}, 120(3):281--295.

\bibitem[Kabluchko et~al., 2009]{kabluchko2009stationary}
Kabluchko, Z., Schlather, M., and de~Haan, L. (2009).
\newblock Stationary max-stable fields associated to negative definite
  functions.
\newblock {\em The Annals of Probability}, 37(5):2042--2065.

\bibitem[McNeil et~al., 2015]{QRM2015}
McNeil, A.~J., Frey, R., and Embrechts, P. (2015).
\newblock {\em Quantitative Risk Management: Concepts, Techniques and Tools}.
\newblock Princeton University Press.

\bibitem[Moltchanov, 2012]{moltchanov2012distance}
Moltchanov, D. (2012).
\newblock Distance distributions in random networks.
\newblock {\em Ad Hoc Networks}, 10(6):1146--1166.

\bibitem[Padoan et~al., 2010]{padoan2010likelihood}
Padoan, S.~A., Ribatet, M., and Sisson, S.~A. (2010).
\newblock Likelihood-based inference for max-stable processes.
\newblock {\em Journal of the American Statistical Association},
  105(489):263--277.

\bibitem[Resnick, 1987]{resnickextreme}
Resnick, S.~I. (1987).
\newblock {\em Extreme Values, Regular Variation, and Point Processes}.
\newblock Springer.

\bibitem[Schlather, 2002]{schlather2002models}
Schlather, M. (2002).
\newblock Models for stationary max-stable random fields.
\newblock {\em Extremes}, 5(1):33--44.

\bibitem[Schlather and Tawn, 2003]{schlather2003dependence}
Schlather, M. and Tawn, J.~A. (2003).
\newblock A dependence measure for multivariate and spatial extreme values:
  Properties and inference.
\newblock {\em Biometrika}, 90(1):139--156.

\bibitem[Serfling, 2002]{serfling2002quantile}
Serfling, R. (2002).
\newblock Quantile functions for multivariate analysis: approaches and
  applications.
\newblock {\em Statistica Neerlandica}, 56(2):214--232.

\bibitem[Smith, 1990]{smith1990max}
Smith, R.~L. (1990).
\newblock Max-stable processes and spatial extremes.
\newblock {\em Unpublished manuscript, University of North Carolina}.

\bibitem[Spodarev, 2014]{spodarev2014limit}
Spodarev, E. (2014).
\newblock Limit theorems for excursion sets of stationary random fields.
\newblock In {\em Modern Stochastics and Applications}, pages 221--241.
  Springer.

\bibitem[Strokorb et~al., 2015]{strokorb2015tail}
Strokorb, K., Ballani, F., and Schlather, M. (2015).
\newblock Tail correlation functions of max-stable processes.
\newblock {\em Extremes}, 18(2):241--271.

\bibitem[Swanson, 2011]{swanson2011band}
Swanson, D. (2011).
\newblock The band around a convex body.
\newblock {\em The College Mathematics Journal}, 42(1):15--24.

\bibitem[Yuen and Stoev, 2013]{yuen2013crps}
Yuen, R. and Stoev, S. (2013).
\newblock Crps {M}-estimation for max-stable models.
\newblock {\em Extremes}, 17(3):387--410.

\end{thebibliography}

\end{document}